\documentclass{article}
\usepackage{amsmath,amsthm,amssymb,amsfonts,mathtools,bm,bbm,array,float,mathdots,dsfont,mathrsfs}
\usepackage{caption,subcaption,multirow,longtable,lscape,wrapfig,tikz,comment}
\usepackage{hyperref}
\hypersetup{allcolors=black}
\usepackage{graphicx,tikz,tikz-cd}
\usetikzlibrary{positioning}
\tikzset{
  node/.style={circle, draw=blue!70!black, fill=blue!20, thick, minimum size=4.5mm, inner sep=0pt},
  edge/.style={->, line width=0.9pt, >=Latex, shorten >=1pt, shorten <=1pt},
  every label/.style={font=\small}
}
\usepackage[top=1in, bottom=1in, left=1in, right=1in]{geometry}
\usepackage[absolute]{textpos}
\usepackage[toc,page]{appendix}
\usepackage[utf8]{inputenc}

\allowdisplaybreaks

\begin{document}
\title{\bf Minimising Willmore Energy via Neural Flow}
\author{\bf Edward Hirst\footnote{ehirst@unicamp.br}, Henrique N. S\'a Earp\footnote{hqsaearp@unicamp.br}, Tom\'as S. R. Silva\footnote{tomas@ime.unicamp.br} }
\date{\small\today\\ \mbox{}\\
\textit{\small Instituto de Matemática, Estatística e Computação Científica (IMECC), \\ Universidade Estadual de Campinas (UNICAMP),\\ 13083-859, Brazil}\\
}

\maketitle

\begin{abstract}
The neural Willmore flow of a closed oriented $2$-surface in $\mathbb{R}^3$ is introduced as a natural evolution process to minimise the Willmore energy, which is the squared $L^2$-norm of mean curvature. Neural architectures are used to model maps from topological $2d$ domains to $3d$ Euclidean space, where the learning process minimises a PINN-style loss for the Willmore energy as a functional on the embedding. Training reproduces the expected round sphere for genus $0$ surfaces, and the Clifford torus for genus $1$ surfaces, respectively. Furthermore, the experiment in the genus $2$ case provides a novel approach to search for minimal Willmore surfaces in this open problem.
\end{abstract}

\thispagestyle{empty}
\clearpage
\setcounter{page}{1}
\numberwithin{equation}{section}
\numberwithin{figure}{section}
\numberwithin{table}{section}
\setcounter{tocdepth}{2}
\tableofcontents
\clearpage 

\section{Introduction}

The Willmore energy of a closed embedded surface in $\mathbb{R}^3$ is a natural quantity, which is directly related to surface tension in many practical applications. 
Since its introduction by T.J. Willmore \cite{Willmore1965}, the analytic study of this energy as a functional on the embedding has been a central topic in differential geometry, not least in the last decade, since the celebrated proof by Marques and Neves \cite{MarquesNeves2014} that the genus $1$ energy is indeed minimised by the Clifford torus. However, their sophisticated analytic and topological techniques are not immediately applicable to higher genus cases, which remain conjectural to this day and whose rigorous solution will likely require new theoretical developments in geometric analysis.

Across the epistemological landscape, AI-driven  techniques have produced  breakthroughs throughout the sciences, albeit so far with comparably limited impact in mathematics. At their core, AI methods — especially supervised ones — require vast amounts of data and operate by minimising a loss function through efficient differentiation. Curvature is a fundamentally differential quantity, making its study particularly amenable to such methods; moreover, the continuous nature of manifolds provides an essentially inexhaustible supply of data that can be sampled cheaply. And yet, the differential geometry community has barely begun to explore the potential of AI techniques as conjectural support devices, let alone as proof assistants. This paper aims to convince that such intersectionality is not only well-motivated but technically and culturally accessible both to geometers and to the AI-assisted maths community. 

In order to adapt a supervised architecture to a differential geometry problem, a reliable procedure is to set up a Physics-Informed Neural Network (PINN) loss, where the functional being minimised is changed from a traditional supervised loss to one inspired by the intrinsic geometry of the problem.
PINNs have already seen an array of early successes on such instances, most notably metric learning, including examples for Einstein metrics on spheres \cite{Hirst:2025seh, Cortes:2026kfx}, Calabi--Yau manifolds \cite{douglas2020numerical, Larfors:2021pbb, Gerdes:2022nzr, Berglund:2022gvm}, and $\mathrm{G}_2$-manifolds \cite{douglas2024harmonic, Heyes:2026rch}; with examples where AI work led to analytically exact computer-assisted proofs \cite{GomezSerranoSurvey, wang2023, platt_nirenberg}. 
Broader applications are still quite incipient, with some impactful first looks at higher-dimensional minimal surfaces \cite{ZhouYe2023MinimalSurfacePINN}, or minimal surfaces in spacetime \cite{Hashimoto:2025zmi}.

\begin{figure}[H]
    \centering
    \begin{subfigure}[t]{0.32\textwidth}
        \centering
        \includegraphics[width=\linewidth]{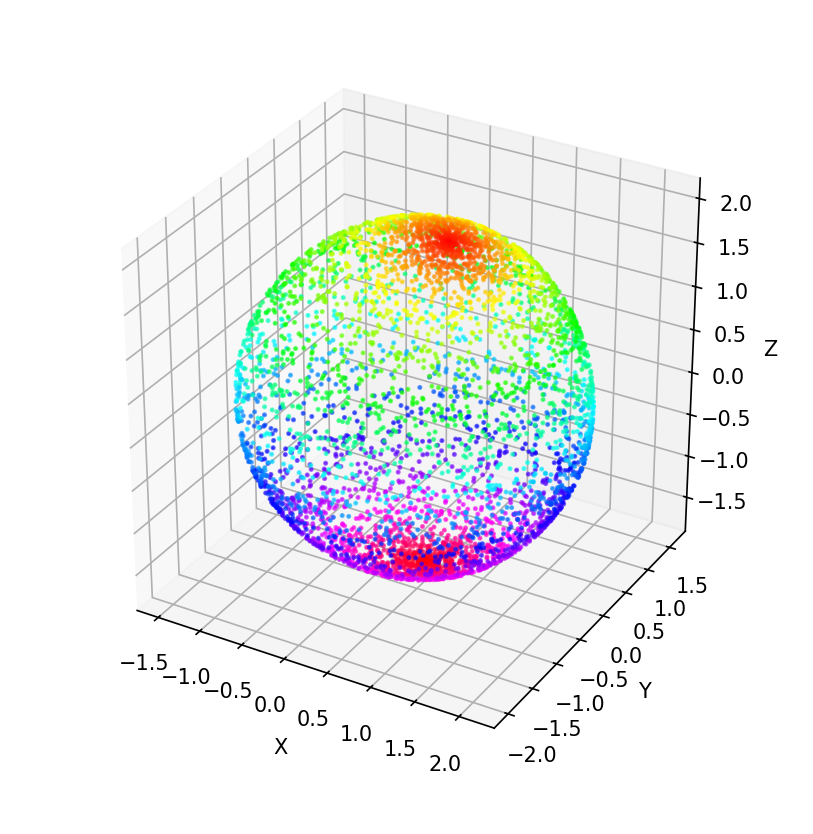}
        \caption{Genus 0}
        \label{fig:g0_final}
    \end{subfigure}
    \hfill
    \begin{subfigure}[t]{0.32\textwidth}
        \centering
        \includegraphics[width=\linewidth]{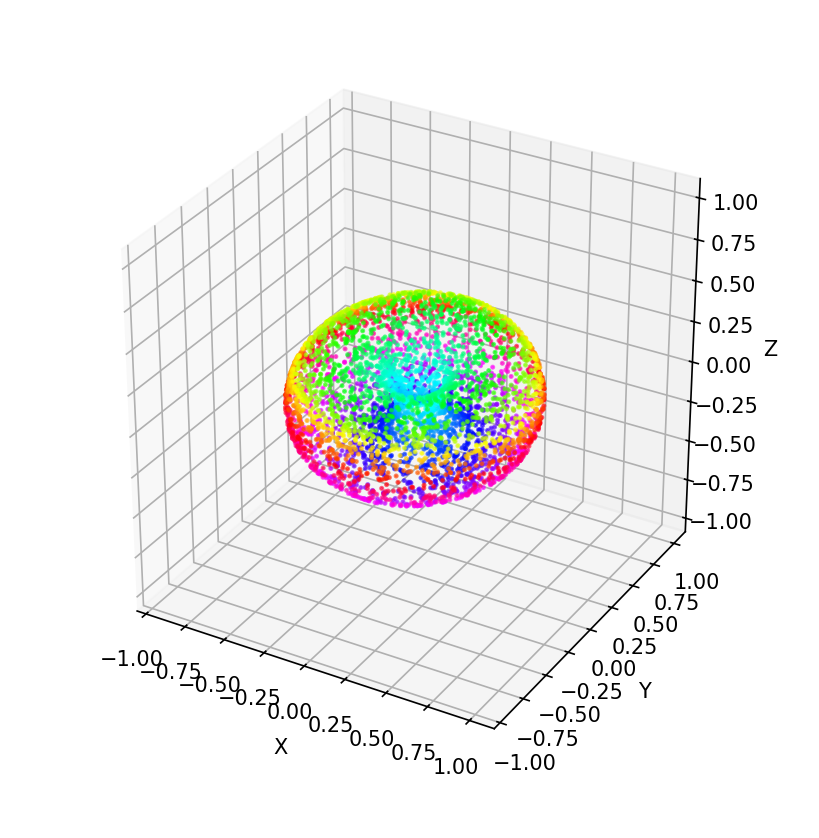}
        \caption{Genus 1}
        \label{fig:g1_final}
    \end{subfigure}
    \hfill
    \begin{subfigure}[t]{0.32\textwidth}
        \centering
        \includegraphics[width=\linewidth]{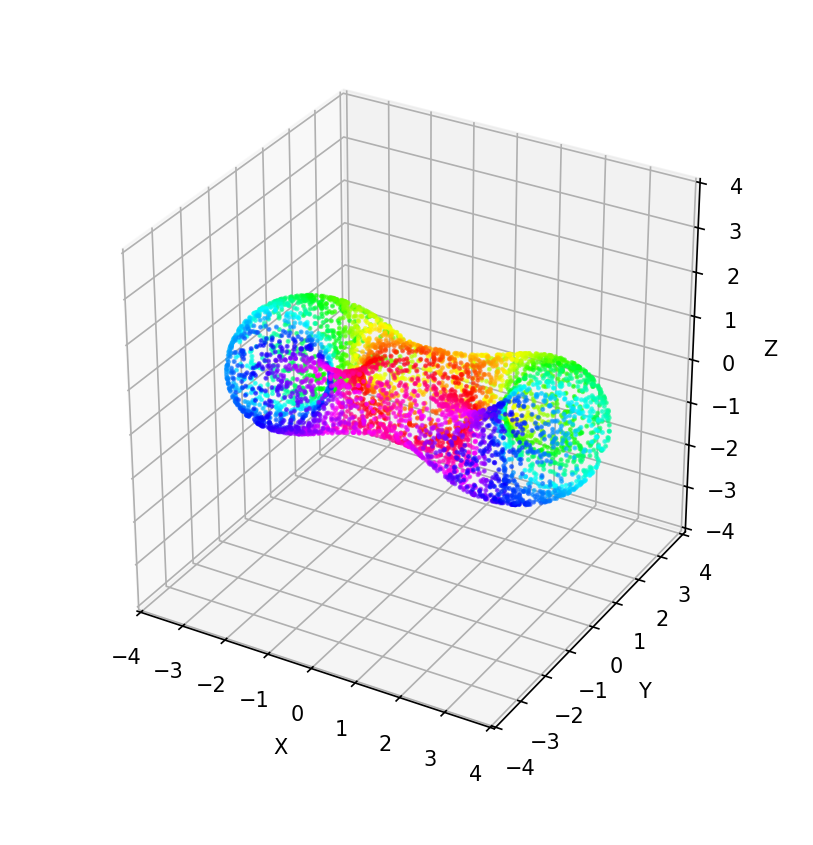}
        \caption{Genus 2}
        \label{fig:g2_final}
    \end{subfigure}
    \caption{The final trained embedded surfaces for each considered genus $\{0, 1, 2\}$, from the neural Willmore flow learning.}
    \label{fig:final_surface}
\end{figure}

Motivated by these successes of PINN-inspired neural architectures, this work investigates the first application of these methods to mean curvature problems via the Willmore energy functional.
The architecture learns an embedding map into $\mathbb{R}^3$ from 2d fundamental domains, for surfaces of genus $\{0, 1, 2\}$, using the highly-optimised \texttt{autograd} routines to efficiently calculate curvatures in a mesh-free manner over millions of points on the surfaces.
Then optimising the embedding to minimise this Willmore energy functional through the machine learning process. Final visualisations of the learnt surfaces achieved in this work, for genuses $\{0,1,2\}$, are shown in Figure \ref{fig:final_surface}.

Computational work is completed in \texttt{python}, using \texttt{pytorch} \cite{torch}, and scripts are made available at this work's respective \href{https://github.com/edhirst/WillmorePINN}{\texttt{GitHub}}\footnote{\href{https://github.com/edhirst/WillmorePINN}{\texttt{https://github.com/edhirst/WillmorePINN}}} repository. For a hands-on introduction, an \emph{interactive tutorial} is available as a \href{https://mybinder.org/v2/gh/edhirst/WillmorePINN/HEAD?urlpath=%2Fdoc%2Ftree%2Fdemo.ipynb}{Binder notebook}.

\section{Background}\label{sec:bkg}
In this section the mathematics behind the Willmore energy is introduced, with motivation and known results.
Then followed by introduction of the relevant AI tools used in this work, specifically PINNs, implemented to numerically model the Willmore energy minimisation.

\subsection{Willmore energy}\label{sec:willmore_bkg}
Let $\Sigma$ be a smooth, closed (compact without boundary), orientable surface of genus $g \geq 0$ in Euclidean space.
An \emph{immersion} of that surface is a smooth map $\varphi \colon \Sigma \to \mathbb{R}^3$ whose differential is everywhere injective, and an \emph{embedding} is a proper injective  immersion (in particular, without self-intersections).
Working in local coordinates $(u, v)$ on a chart $U \subset \mathbb{R}^2$, we write $\varphi = \varphi(u,v)$ and denote partial derivatives by $\varphi_u = \partial_u\varphi$ and $\varphi_v = \partial_v\varphi$.

The \emph{first fundamental form} (induced metric) is the symmetric covariant $2$-tensor
\begin{equation}
    \mathrm{I} = E\,du^2 + 2F\,du\,dv + G\,dv^2,
\end{equation}
with coefficients
\begin{equation}
    E = \langle \varphi_u,\, \varphi_u \rangle, \qquad
    F = \langle \varphi_u,\, \varphi_v \rangle, \qquad
    G = \langle \varphi_v,\, \varphi_v \rangle,
\end{equation}
where $\langle \cdot, \cdot \rangle$ denotes the standard Euclidean inner product on $\mathbb{R}^3$.
The immersion condition requires $$\det \mathrm{I} = EG - F^2 > 0$$ everywhere, and it defines the induced area element as $dA = \sqrt{EG - F^2}\,du\,dv$.  
The unit normal field is then
\begin{equation}
    \hat{n} = \frac{\varphi_u \times \varphi_v}{\lvert \varphi_u \times \varphi_v \rvert}.
\end{equation}
The \emph{second fundamental form} then encodes the extrinsic curvature of the immersion:
\begin{equation}
  \mathrm{II} = L\,du^2 + 2M\,du\,dv + N\,dv^2,
\end{equation}
where
\begin{equation}
    L = \langle \varphi_{uu},\, \hat{n} \rangle, \qquad
    M = \langle \varphi_{uv},\, \hat{n} \rangle, \qquad
    N = \langle \varphi_{vv},\, \hat{n} \rangle.
\end{equation}
These forms are used to define the shape operator $\mathcal{S} \colon T_p\Sigma \to T_p\Sigma$, which is then given by $\mathcal{S} = \mathrm{I}^{-1}\mathrm{II}$ at each point $p \in \Sigma$.
Its eigenvalues $\kappa_1, \kappa_2$ are the principal curvatures, which relate to the two most standard curvature measures: the \emph{mean curvature}  $H$ and the \emph{Gaussian curvature} $K$, computable in the chart $\varphi$ by
\begin{equation}
    H = \frac{\kappa_1 + \kappa_2}{2} = \frac{EN - 2FM + GL}{2(EG - F^2)},
    \qquad\qquad
    K = \kappa_1\kappa_2 = \frac{LN - M^2}{EG - F^2}.
\end{equation}

In this article we will be solely interested in the mean curvature, since indeed the \emph{Willmore energy} of an immersed surface $\varphi \colon \Sigma \to \mathbb{R}^3$ is defined by
\begin{equation}
  \mathcal{W}(\varphi) = \int_\Sigma H^2\,dA.
  \label{eq:willmore}
\end{equation}
This functional was systematically studied by Willmore \cite{Willmore1965}, with precursors in the nineteenth-century elasticity literature and in the work of Thomsen \cite{Thomsen1923}.
It measures, in an $L^2$ sense, the deviation of the surface from a minimal surface ($H \equiv 0$).
The Willmore energy \eqref{eq:willmore} is scale-invariant, under a uniform rescaling $\varphi \mapsto \lambda\varphi$ for $\lambda > 0$, the mean curvature transforms as $H \mapsto \lambda^{-1} H$ while the area element transforms as $dA \mapsto \lambda^2\,dA$, so that
\begin{equation}
    \mathcal{W}(\lambda\varphi)
    = \int_\Sigma (\lambda^{-1}H)^2 (\lambda^2\,dA)
    = \int_\Sigma H^2\,dA
    = \mathcal{W}(\varphi).
\end{equation}
This is a special case of the full Möbius invariance $S^3 = \mathbb{R}^3 \cup \{\infty\}$ \cite{Blaschke1929, White1973}: uniform scaling is a conformal diffeomorphism of $\mathbb{R}^3$, and the Willmore energy is invariant under the entire conformal group.  
Crucially, this invariance holds for immersed surfaces of \emph{any} genus: neither the energy value nor the location of its critical points depends on the global scale of the embedding.  
In particular, the known minimum values $4\pi$ (genus 0) and $2\pi^2$ (genus 1) are scale-free constants.

As a side note, since the Gauss--Bonnet theorem gives $\int_\Sigma K\,dA = 2\pi\chi(\Sigma)$, which is a topological invariant, minimising $\mathcal{W}$ over immersions of a fixed genus $g$ is also equivalent to minimising the \emph{conformal Willmore energy}
\begin{equation}
    \widetilde{\mathcal{W}}(\varphi)
    = \int_\Sigma (H^2 - K)\,dA
    = \mathcal{W}(\varphi) - 2\pi\chi(\Sigma).
\end{equation}
Immersions that are critical points of $\mathcal{W}$ satisfy the \emph{Willmore equation}, the Euler--Lagrange condition
\begin{equation}
  \Delta_g H + 2H(H^2 - K) = 0,
  \label{eq:willmore-eq}
\end{equation}
where $\Delta_g$ denotes the Laplace--Beltrami operator of the induced metric; such surfaces are called \emph{Willmore surfaces}. 
Traditionally, solving this Euler--Lagrange equation is the approach modelling Willmore flow analytically, and it is this step that this work replaces with a machine learning method, directly minimising $\mathcal{W}$ as a neural flow \cite{Halverson:2023ndu}.

From a computational perspective, the Willmore functional acts as a curvature-regularising prior, it penalises high-frequency geometric oscillations while remaining scale-invariant and invariant under conformal reparametrisation.  
In contrast to area minimisation, which risks topological collapse, Willmore-type energies control curvature concentration and promote smoothness without inducing shrinkage.  
These properties make the functional a principled geometric objective for learning surface embeddings via neural architectures, where robustness under rescaling and reparametrisation is essential.  
The present work exploits this structure directly: an immersion $\varphi \colon \Sigma_g \to \mathbb{R}^3$ is learnt by minimising $\mathcal{W}$ through automatic differentiation of the fundamental forms, with topology fixed by construction via the parametric domain of genus $g$.

\subsubsection{Known minimisers by genus}
The problem of minimising $\mathcal{W}$ over immersed surfaces of prescribed genus $g$ has been resolved for $g \leq 1$, with the higher-genus case remaining largely conjectural.

\paragraph{Genus $0$.}
For immersed spheres, Willmore \cite{Willmore1965} established that $\mathcal{W}(\varphi) \geq 4\pi$ for every smooth immersion $\varphi \colon S^2 \to \mathbb{R}^3$, with equality if and only if $\varphi$ parametrises a round sphere, up to Möbius transformation.  
This lower bound also follows from the Li--Yau inequality \cite{LiYau1982}.  
The existence of a smooth energy-minimising sphere was proved by Simon \cite{Simon1993}.

\paragraph{Genus $1$.}
For tori, Willmore conjectured \cite{Willmore1965} that
\begin{equation}
    \mathcal{W}(\varphi) \geq 2\pi^2
\end{equation}
for every smooth immersion $\varphi \colon T^2 \to \mathbb{R}^3$.  The conjectured minimiser is the \emph{Clifford torus}, the product of two circles of radius $1/\sqrt{2}$ embedded in $S^3 \subset \mathbb{R}^4$ and stereographically projected to $\mathbb{R}^3$, equivalently the standard torus of revolution with major and minor radii satisfying $R/r = \sqrt{2}$, for which $\mathcal{W} = 2\pi^2$. 
After standing open for nearly fifty years, the conjecture was proved in full generality by Marques and Neves \cite{MarquesNeves2014} using min-max theory in the space of flat varifolds.

\paragraph{Genus $g \geq 2$.}
For surfaces of higher genus, Bauer and Kuwert \cite{BauerKuwert2003} established the existence of smooth Willmore minimisers; see also \cite{Simon1993, Riviere2008} for related compactness and regularity results.   
Knowing that the Clifford torus is the global Willmore minimiser across genuses, and that catenoid bridges have zero Willmore energy, for each genus one can build a surface of genus $g>1$ by gluing $g$ Clifford tori via Catenoid bridges. 
The Willmore energy via this construction would then be $g(2\pi^2)$, which for genus 2 would give $4\pi^2 \sim 39.48$.

The precise infimum of $\mathcal{W}$ over genus-$g$ immersions remains conjectural for $g \geq 2$.  
Kusner \cite{Kusner1989} conjectured that the minimisers are the Lawson minimal surfaces $\xi_{g,1} \subset S^3$ of \cite{Lawson1970}, stereographically projected to $\mathbb{R}^3$. 
The conjectured minimiser $\xi_{2,1}$ has numerically computed Willmore energy \cite{Hsu1992}
\begin{equation}
    \inf_{g=2}\,\mathcal{W} = \mathcal{W}(\xi_{2,1}) \approx 4\pi^2 \approx 21.89.
\end{equation}

\subsection{PINNs}\label{sec:pinns_bkg}
The reader unfamiliar with machine learning may wish to consult Appendix~\ref{appendix:ml} for a gentle introduction before proceeding.

A neural network is a parametric function $f_\theta \colon \mathbb{R}^{d_{\mathrm{in}}} \to
\mathbb{R}^{d_{\mathrm{out}}}$ formed by composing affine maps with pointwise nonlinear
activation functions $\sigma$ (e.g.\ $\tanh$, ReLU):
\begin{equation}
    f_\theta(x) = W_L\,\sigma\!\bigl(W_{L-1}\,\sigma(\cdots\sigma(W_1 x + b_1)\cdots) + b_{L-1}\bigr) + b_L,
\end{equation}
where the weights $\theta = \{W_\ell, b_\ell\}$ are the trainable parameters \cite{Goodfellow2016}.
The architecture hyperparameters (network depth $L$, layer widths, and choice of $\sigma$) are fixed before training and determine the expressiveness of the function class.
Given a \emph{loss function} $\mathcal{L}(\theta)$, traditionally measuring discrepancy between the network output and a target objective, parameters are optimised by stochastic gradient descent: at each iteration a random mini-batch of $B$ input points is drawn, the loss and its gradient $\nabla_\theta \mathcal{L}$ are computed via backpropagation \cite{Rumelhart1986}, and $\theta$ is updated by an optimiser (here Adam \cite{Kingma2015}) that maintains per-parameter adaptive learning rates.
Iterating over the full dataset defines one epoch; performance on a held-out test set monitors generalisation.

Importantly, a feedforward neural network with a single hidden layer and a non-polynomial activation function can approximate any continuous function on a compact domain to arbitrary precision \cite{Cybenko1989, Hornik1991}.  
Beyond this tautological case, deeper architectures with fixed width extend this result and offer substantially more favourable approximation rates for smooth functions \cite{Barron1993, DeVore2021}.  
More precisely, for a target function $f \in C^k(\Omega)$, networks with $L$ layers of width $N$ achieve approximation error $O(N^{-k/d})$ in $d$ dimensions \cite{Yarotsky2017,
DeVore2021}, recovering classical polynomial rates while exploiting architecture structure to partially alleviate the curse of dimensionality.

These approximation-theoretic guarantees from neural architectures motivate representing unknowns in variational problems as the output of a neural network, and replacing the variational criterion with a differentiable training objective.  
This paradigm, broadly termed \emph{physics-informed neural networks} (PINNs), was introduced in \cite{Raissi2019} for PDE-constrained problems, and has since been applied to a wide class of variational and inverse problems \cite{Karniadakis2021}.  
The key observation is that automatic differentiation \cite{Baydin2018} allows partial derivatives of the network output with respect to its inputs to be computed efficiently and exactly, to machine precision.
Allowing differential operators appearing in physical or geometric equations to be evaluated and differentiated through without discretisation error.
It is the use of these derivative terms in loss functions which transforms training of a NN architecture into training of a PINN architecture.

\subsubsection{PINNs for geometric variational problems}
The Willmore energy \eqref{eq:willmore} is a second-order\footnote{Whereas its Euler--Lagrange equation \eqref{eq:willmore-eq} is fourth-order.} geometric functional, involving second derivatives of the immersion $\varphi$.  
Classical numerical approaches — finite elements on triangulated meshes, or discrete differential geometry formulations — must either enforce $C^1$ continuity across element boundaries or resort to mixed formulations to handle the biharmonic structure \cite{Dziuk2008, Crane2013}.  
Moreover, mesh-based methods fix a discrete combinatorial
structure that constrains the topology and resolution of the solution.

A PINN formulation bypasses these difficulties naturally.  
The immersion
$\varphi \colon \Omega \to \mathbb{R}^3$ is represented as a smooth function
\begin{equation}
    \varphi_\theta \colon \Omega \subset \mathbb{R}^2 \longrightarrow \mathbb{R}^3,
\end{equation}
parametrised by network weights $\theta \in \mathbb{R}^p$.  
The first and second fundamental forms are computed pointwise via automatic differentiation of $\varphi_\theta$, yielding exact expressions for the coefficients $E, F, G, L, M, N$ and hence for $H$ and $dA$ at any sample point in $\Omega$.  
The training objective is then a Monte Carlo approximation of $\mathcal{W}(\varphi_\theta)$, and gradient descent in $\theta$ performs continuous gradient flow in the space of immersions represented by the network class.

This approach offers several structural advantages for the Willmore problem.  
First, the network is intrinsically smooth, so no additional regularisation is required to maintain differentiability of the immersion.  
Second, the parametric domain $\Omega$ is fixed and encodes the topology, since choosing $\Omega$ to be a fundamental domain of genus $g$ with appropriate boundary identifications, the topology of the immersed surface is preserved throughout training without any explicit topological constraint in the loss.  
Third, the conformal invariance of $\mathcal{W}$ is respected at the level of the continuous functional, not only approximately at a discrete level.  
Finally, the method scales naturally to higher-resolution solutions by increasing network width or depth, independently of any mesh refinement strategy.

\section{The neural architecture}\label{sec:architecture}
In this section, details related to the implemented methods for this work are given.
This includes information about sampling the fundamental domains to form the inputs, and how and what losses are computed on the outputs.
There is then discussion of the supervised pretraining used to initiate the architecture to the correct topology in each case before training; followed by the specific architecture details and hyperparameters used for the reported runs.

\subsection{Sampling}
In this work, immersions of surfaces of varying genus are considered.
For each target genus, the neural architecture learns a map $\varphi_\theta \colon \Omega \to \mathbb{R}^3$ for a fixed domain $\Omega$, which reflects the surface's topology.

\paragraph{Genus 0.}
The parameter domain is $\Omega_0 = [0, 2\pi] \times [0, \pi]$, with $u \in [0, 2\pi]$ the azimuthal angle and $v \in [0, \pi]$ the polar angle.
Points are sampled uniformly in $(u, v)$; the area element $\sqrt{EG-F^2}$ appearing explicitly in the Monte Carlo integrand accounts for the induced metric without importance reweighting.

To enforce the topological constraint that both poles collapse to single points, the network input uses spectral features (inspired by \cite{Berglund:2022gvm}), which here are real spherical harmonics.
For degree $l$ and order $m$, the real harmonics are
\begin{equation}
    Y_l^m(\theta, \phi) = c_{l,m}\,P_l^{|m|}(\cos\theta) \cdot
    \begin{cases}
    \sin(|m|\phi) & m < 0, \\
    1             & m = 0, \\
    \cos(|m|\phi)   & m > 0,
    \end{cases}
\end{equation}
where $P_l^m$ are the associated Legendre polynomials and $c_{l,m}$ the standard orthonormality constants.  Pole regularity is automatic: since $P_l^m(\pm 1) = 0$ for $m \neq 0$, all $\phi$-dependent harmonics vanish at $\theta = 0$ and $\theta = \pi$, ensuring the network output is $u$-independent at both poles.  Using harmonics up to degree $L$ yields $(L+1)^2$ input features.

\paragraph{Genus 1.}
The parameter domain is $\Omega_1 = [0, 2\pi]^2$, doubly periodic.
Points are sampled uniformly on the unit square and mapped to a parallelogram fundamental domain via the affine transformation
\begin{equation}
  (u, v) \;\mapsto\;
  \bigl(u + v\,\mathrm{Re}(\tau),\;\; v\,\mathrm{Im}(\tau)\bigr) \cdot 2\pi,
  \qquad \tau \in \mathbb{C},\quad \mathrm{Im}(\tau) > 0.
\end{equation}
Periodicity in both directions is encoded by a Fourier feature embedding: for $k = 1, \ldots, N$,
\begin{equation}
    \Phi(u, v) = \bigl(\sin ku,\,\cos ku,\,\sin kv,\,\cos kv\bigr)_{k=1}^{N},
\end{equation}
giving $4N$ features.  
Since $\sin(k \cdot 2\pi j) = 0$ and $\cos(k \cdot 2\pi j) = 1$ for integer $j$, the feature vector is identical at all identified boundary points, enforcing doubly-periodic boundary conditions exactly in the input representation.

\paragraph{Genus 2.} 
In an ideal scenario higher genus surfaces would be represented by a single connected fundamental domain, with $4g$ edges and appropriate pairwise identifications. 
Inputs would then be the corresponding Laplacian eigenfunctions, which would enforce these identifications intrinsically.
For the genus 2 case this would be an octagon fundamental domain, and although a large amount of code was developed along this line (available in other branches of the repository), a key obstruction was that since there is no known analytic embedding function from this octagon to a genus 2 surface in $\mathbb{R}^3$, the genus 2 nature could not be initialised for the architecture and learning could not progress.
It it worth noting too that the Laplacian eigenfunctions also don't have analytic forms, and had to be numerically approximated.

In place of this, an alternative approach for the genus 2 construction and sampling was taken.
Here the genus 2 surface is represented by two torus charts $T_1, T_2$, each with domain $\Omega = [0,2\pi]^2$, glued along boundary circles.
For this, on each chart a disc of parameter-space radius $\delta = 0.65$, centred at $(0, 0)$ and $(0, \pi)$ respectively, is excluded via rejection sampling, i.e. candidate points $(u,v) \sim \mathrm{Uniform}([0,2\pi]^2)$ are discarded whenever
\begin{equation}
  \min_{n_u, n_v \in \{0,\pm 1\}} \bigl((u - 2\pi n_u)^2 + (v - 2\pi n_v)^2\bigr)
  \;\leq\; \delta^2,
\end{equation}
where the minimum accounts for periodic wrap-around.
In addition to the bulk samples, an annular transition zone is set up at the border of the excluded discs, defined via
$\delta \leq r_{u,v} \leq \alpha\delta$ (with $\alpha = 2.5$)
on each chart, which is sampled separately and uniformly for a gluing loss in the training and an annular regularity term added to the regularity loss.
Doubly-periodic boundary conditions on each chart are enforced by the same Fourier feature embedding as genus 1, again giving $4N$ features per chart.

In defining the gluing map between the annuli surrounding the excluded discs of each torus fundamental domain, each chart carries a local polar coordinate system centred at its disc centre:
$(0, 0)$ on $T_1$ and $(\pi, 0)$ on $T_2$.
For a point $(u,v)$ on chart $T_i$, the signed periodic displacement from the centre is $\tilde{u} = [(u - u_0^{(i)} + \pi) \bmod 2\pi] - \pi$ (and likewise $\tilde{v}$), giving polar coordinates: $(r, \theta) = (\sqrt{\tilde{u}^2 + \tilde{v}^2}, \mathrm{atan2}(\tilde{v},\, \tilde{u}))$.
The two charts are then identified by the radial reflection
\begin{equation}
  \Theta \colon (r, \theta)_{T_1} \;\longmapsto\; (2\delta - r,\, \theta)_{T_2},
\end{equation}
so that matched sample pairs are
\begin{equation}
  p_1 = \bigl(r\cos\theta,\; r\sin\theta\bigr), \qquad
  p_2 = \bigl(\pi + (2\delta - r)\cos\theta,\; (2\delta - r)\sin\theta\bigr),
\end{equation}
with all coordinates taken mod $2\pi$, $r \sim \mathrm{Uniform}[\delta(1-w),\, \delta(1+w)]$, and $\theta \sim \mathrm{Uniform}[0, 2\pi)$.
The Jacobian of $\Theta$ satisfies $D\Theta\colon e_r \mapsto -e_r$ and
$D\Theta\colon e_\theta \mapsto \tfrac{2\delta-r}{r}\,e_\theta$, which later determines the $C^1$ and $C^2$ matching conditions enforced by the gluing loss.

\subsection{Loss functions}
During the learning process the immersion map, defined by the neural architecture and parameterised by the model parameters $\theta$, updates as the model parameters are changed and the training loss is minimised.
This training loss is a weighted combination of a Willmore energy term, a composite regularity term, and for genus 2 a gluing term also.
Set up is such that the learning process models a flow in the Willmore energy whilst maintaining regularity conditions for a consistent embedding of the surface.

\paragraph{Willmore Loss.}
The Willmore functional is approximated by Monte Carlo integration over a batch of $B$ uniform parameter samples $\{(u_i, v_i)\}_{i=1}^B$:
\begin{equation}
    \mathcal{L}_\mathcal{W}(\varphi_\theta)
    = |\Omega|\,\frac{1}{B}\sum_{i=1}^{B}
    H_i^2 \,\sqrt{E_i G_i - F_i^2},
    \label{eq:mc-willmore}
\end{equation}
where $|\Omega|$ is the area of the parameter domain, $E_i, F_i, G_i$ and $L_i, M_i, N_i$ are the fundamental form coefficients, and $H_i$ the full mean curvature, at $(u_i, v_i)$, computed by automatic differentiation of the model $\varphi_\theta$.
In this and all following losses L2 norms over the output $\mathbb{R}^3$ components are also assumed where appropriate.

For genus-2 surfaces, the two-chart architecture parametrises each punctured torus $T_k \setminus D_k$ ($k=1,2$) independently, where $D_k$ is the small excluded disc of parameter-space radius $\delta$ that is removed to accommodate the gluing neck.
The Willmore estimate \eqref{eq:mc-willmore} is applied separately to each chart with domain area $|\Omega_k| = (2\pi)^2 - \pi\delta^2$, and the contributions are summed, $\mathcal{L}_\mathcal{W} = \mathcal{L}_{\mathcal{W}_{T_1}} + \mathcal{L}_{\mathcal{W}_{T_2}}$.
To stabilise training, a Huber-type smoothing is applied to the $H^2$ integrand before multiplying by the area element.
For a threshold $c$ (set so that values exceeding $\sim 50\times$ the batch mean are clipped in practice), the modified integrand is
\begin{equation}
    \tilde{H}^2 =
    \begin{cases}
        H^2, & H^2 \leq c, \\
        2\sqrt{c}\,|H| - c, & H^2 > c,
    \end{cases}
\end{equation}
which matches the exact Willmore integrand—and hence its minimiser—for bounded $H$, but replaces the quadratic growth above the threshold with a linear ramp, keeping gradients constant in magnitude rather than vanishing at high-curvature junctions.
In practice the threshold was only exceeded for early stages of training where the gluing junction was not smooth, and the standard $\mathcal{L}_{\mathcal{W}}$ (without this Huber modification) was computed separately for all reported losses and testing.


\paragraph{Regularity Loss.}
To prevent degeneracy during optimisation, a regularity loss $\mathcal{L}_R$ penalises three failure modes.  
Let $a_{\min}$ be a prescribed minimum area element, $\epsilon_{\text{pos}}$ a small lower threshold on the diagonal metric components, and $g_{\max}$ a maximum metric value.
The three components are:
\begin{align}
  \mathcal{L}_{\text{area}}
    &= \frac{1}{B}\sum_i \bigl[a_{\min}
       - \sqrt{E_i G_i - F_i^2}\bigr]_+^2,  \label{eq:reg-area} \\
  \mathcal{L}_{\text{pos}}
    &= \frac{1}{B}\sum_i \bigl([\epsilon_{\text{pos}} - E_i]_+^2 + [\epsilon_{\text{pos}} - G_i]_+^2\bigr),\\
  \mathcal{L}_{\text{smooth}}
    &= \frac{1}{B}\sum_i \bigl([E_i - g_{\max}]_+^2 + [G_i - g_{\max}]_+^2\bigr),
\end{align}
where $[x]_+ = \max(0, x)$ is a ReLU gate, and means if regularity violation is negligible then these losses stay at 0 and do not contribute, avoiding interference with Willmore minimisation when the metric is well-conditioned.  
The regularity loss is the normalised sum
\begin{equation}
  \mathcal{L}_R = \frac{\alpha_{\text{area}} \mathcal{L}_{\text{area}}
                      + \alpha_{\text{area}} \mathcal{L}_{\text{pos}}
                      + \alpha_{\text{smooth}} \mathcal{L}_{\text{smooth}}}
                       {\alpha_{\text{area}} + \alpha_{\text{pos}} + \alpha_{\text{smooth}}}.
\end{equation}

For the genus 2 training the regularity loss is computed on each fundamental domain chart and summed, with two additional modifications.
First, since the Willmore and regularity integrals exclude the disc $D_k$ of parameter-space radius $\delta$ on each chart $T_k$, the network is otherwise free to form an arbitrarily thin funnel just outside the gluing circle.
To prevent this, $\mathcal{L}_R$ is evaluated on an additional annular sample $\delta \leq r_k \leq \alpha\delta$, and a weighted copy is added to the total regularity.
Second, a log-barrier term is added to $\mathcal{L}_R$:
\begin{equation}
    \mathcal{L}_{\mathrm{log}}
    = \frac{1}{B}\sum_i
      \Bigl[\,{-\log\sqrt{E_i G_i - F_i^2}}\Bigr]_+,
\end{equation}
which acts similarly to $\mathcal{L}_{\text{area}}$, but with $[\,\cdot\,]_+$ restricting the barrier to area elements below 1 (the collapse direction).  
Unlike the previous ReLU-gated terms, whose gradients are zero inside the safe band, this term has a continuous, everywhere non-zero gradient $-1/\!\sqrt{EG-F^2}$ that diverges as the area element approaches zero, providing an always-on signal that forestalls collapse before $\mathcal{L}_{\mathrm{area}}$ begins to fire.

\paragraph{Gluing Loss.}
Since the fundamental domain is disconnected for the genus 2 training, this implementation required an additional loss term, the gluing loss.
The trained architecture must produce a smooth surface across the gluing junction between the two torus charts, and since the Willmore integrand depends on the immersion map only up to second derivatives (via the mean curvature), $C^2$ continuity is precisely the minimal regularity required for the energy to be well-defined and finite on the joined surface.

To enforce this continuous gluing, sampled pairs which represent the same global point in each chart's coordinates (within the gluing annuli) are passed through the embedding network.
Then Jacobians $J_i = D\varphi_i \in \mathbb{R}^{3 \times 2}$ and Hessians $\mathcal{H}_i \in \mathbb{R}^{3 \times 2 \times 2}$ at each sample pair are computed by automatic differentiation through the respective chart sub-networks.
Setting $r_2 = 2\delta - r$ and writing $e_r = (\cos\theta, \sin\theta)^\top$,
$e_\theta = (-\sin\theta, \cos\theta)^\top$ for the polar frame, the three
matching conditions imposed by the gluing map $\Theta$ are:
\begin{align}
  \mathcal{L}_{C^0} &= \frac{1}{B}\sum_i \big(\varphi_1(p_1) - \varphi_2(p_2)\big)^2 , \\
  \mathcal{L}_{C^1} &= \frac{1}{B}\sum_i \big( J_1 e_r + J_2 e_r \big)^2 + \big(r(J_1 e_\theta) - r_2(J_2 e_\theta)\big)^2, \\
  \mathcal{L}_{C^2} &= \frac{1}{B}\sum_i \big(e_r^\top\mathcal{H}_1 e_r - e_r^\top\mathcal{H}_2 e_r\big)^2 \notag \\
                    &\quad + \big(r(e_r^\top\mathcal{H}_1 e_\theta) + r_2(e_r^\top\mathcal{H}_2 e_\theta) + (J_1 e_\theta + J_2 e_\theta)\big)^2 \notag \\
                    &\quad +  \big(r^2(e_\theta^\top\mathcal{H}_1 e_\theta) - r_2^2(e_\theta^\top\mathcal{H}_2 e_\theta) - 2\delta(J_1 e_r)\big)^2 ,
\end{align}
which then combine to form the full gluing loss as
\begin{equation}
  \mathcal{L}_G = \frac{\alpha_{C^0}\,\mathcal{L}_{C^0} + \alpha_{C^1}\,\mathcal{L}_{C^1} + \alpha_{C^2}\,\mathcal{L}_{C^2}}{\alpha_{C^0} + \alpha_{C^1} + \alpha_{C^2}}\;.
\end{equation}

The $C^2$ conditions include curvature-of-frame corrections from $\partial_\theta e_r = e_\theta$ and $\partial_\theta e_\theta = -e_r$, arising when differentiating the $C^0$ identity twice in $(r,\theta)$.
At $r = \delta$ all conditions reduce to the naive $(Dg)^{\otimes k}$ forms.
The total gluing loss is the mean squared residual of each condition, summed with weights $\lambda_{C^0} = 200$, $\lambda_{C^1} = 50$, $\lambda_{C^2} = 20$;.
The $C^1$ and $C^2$ terms are activated with epoch delays to allow $C^0$ continuity, then $C^1$ continuity, to be established sequentially, before entering the full training regime.

\paragraph{Total Loss.}
Combining these loss terms, the total loss is
\begin{equation}
  \mathcal{L} = \frac{\lambda_{\mathcal{W}}(t)\,\mathcal{L}_{\mathcal{W}} + \lambda_R(t)\,\mathcal{L}_R + \lambda_G(t)\,\mathcal{L}_G }
                     {\lambda_{\mathcal{W}}^0 + \lambda_R^0 + \lambda_G^0},
  \label{eq:total-loss}
\end{equation}
normalised by the initial weights $\lambda_{\mathcal{W}}^0, \lambda_R^0, \lambda_G^0$ so the loss scale is stable throughout training (note $\lambda_G = 0$ for genus 0 and 1).  
To adapt the total loss throughout training, allowing fine-tuning towards a good solution, the code has the option for annealing to be applied to the loss weights.
Used only in the genus 2 case, the Willmore weight $\lambda_{\mathcal{W}}(t)$ has an initial warm up period increasing from $\lambda_{\mathcal{W}}^0$ to $\lambda_{\mathcal{W}}^{\text{final}}$, while $\lambda_R(t)$ and $\lambda_G(t)$ decrease equivalently, set up to progressively shift emphasis onto curvature minimisation once a well-conditioned metric has been established.  

\subsection{Supervised Pretraining}
Before the neural Willmore flow is run, the architecture is pretrained to match a closed-form reference embedding.  
This initialises the weights in a topologically correct configuration, reducing the risk of the model converging to a degenerate immersion.

The pretraining loss is different, and combines position and first-derivative matching between the predicted embedding and the reference embedding:
\begin{equation}
  \mathcal{L}_{\text{pre}}(\theta)
  = \frac{1}{B}\sum_{i=1}^B
    \Bigl(
      \beta_1\,\|\varphi_\theta - \varphi_{\text{ref}}\|^2
      + \beta_2(t)\,\bigl(
          \|\partial_u\varphi_\theta - \partial_u\varphi_{\text{ref}}\|^2
        + \|\partial_v\varphi_\theta - \partial_v\varphi_{\text{ref}}\|^2
        \bigr)
    \Bigr),
  \label{eq:pretrain-loss}
\end{equation}
where derivatives of both $\varphi_\theta$ and $\varphi_{\text{ref}}$ are computed by automatic differentiation.  
The code has functionality to allow the derivative weight $\beta_2(t)$ to increase gradually from zero over the first pretraining epochs, allowing the network to first fit position before being constrained on derivatives, although for the presented run results this was not used.

\paragraph{Genus 0 reference.}
The reference surface is the ellipsoid
\begin{equation}
    \varphi_{\text{ref}}(u,v) = (a\sin v\cos u,\; b\sin v\sin u,\; c\cos v), \qquad (u,v) \in [0,2\pi]\times[0,\pi],
\end{equation}
with free parameters $(a, b, c) \in \mathbb{R}_{>0}^3$ prescribing the semi-axes.  
Setting $a = b = c$ recovers a round sphere with Willmore energy exactly $4\pi$.

\paragraph{Genus 1 reference.}
The reference surface is the torus parametrised by the complex modulus $\tau \in \mathbb{C}$, $\mathrm{Im}(\tau) > 0$:
\begin{equation}
  \varphi_{\text{ref}}(u,v)
  = \Bigl(
      (1 + r\cos\tilde v)\cos u,\;
      (1 + r\cos\tilde v)\sin u,\;
      r\sin\tilde v
    \Bigr),
\end{equation}
where $r = |\mathrm{Im}(\tau)|$, the normalised minor-circle coordinate $\tilde v = v/r + (\mathrm{Re}(\tau)/|\mathrm{Im}(\tau)|)\,u$ incorporates a helical twist controlled by $\mathrm{Re}(\tau)$, and the major radius is fixed to unity.  
The free parameter $\tau$ simultaneously determines the minor-to-major radius ratio $r$ and the shear class of the fundamental domain.  
The Clifford torus (minimiser of $\mathcal{W}$ among tori) is approached as $r \to 1/\sqrt{2}$.

\paragraph{Genus 2 reference.}
The reference surface is a disconnected sum of two punctured tori\footnote{Note there is no known analytic expression for an $\mathbb{R}^3$ embedding of the $\xi_{2,1}$ Lawson surface, such that it could not be used as the reference. There is the beginnings of code in a separate branch of the repository which considers the embedding via $S^3$ if one is interested.}.
Each of the $T_k$ ($k=1,2$) chart subnetworks use the same $\varphi_{\text{ref}}(u,v)$ torus reference as in the genus 1 case, parametrised by a complex modulus $\tau_k$ with $\mathrm{Im}(\tau_k)>0$.
The two tori are oriented with their central hole axes along $z$, and displaced along the $x$-axis so that their outer equatorial circles coincide at the origin.
$T_1$ is translated by $-(R+r_1)$ and $T_2$ is translated by $+(R+r_2)$, where $R$ is a fixed centre-shift parameter and $r_k = \mathrm{Im}(\tau_k)$.
This places the facing equators of both charts as the closest points to $x = 0$ with a matching orientation, and the surfaces non-overlapping on either side.
The excluded disc $D_k$ of parameter-space radius $\delta$ is centred at $(u,v)=(0,0)$ for $T_1$ and $(u,v)=(\pi,0)$ for $T_2$ such that the discs face each other and make near contact centred near the origin (in the $y-z$ plane).
Pretraining runs for a longer 200 epochs here to
account for the more complex displaced geometry.

\subsection{Architecture}
The architecture used in this work is a neural architecture in nature, the losses used in training are PINN-inspired, and the architecture uses traditional feedforward neurons at its core.
In the genus 0 and 1 cases the architecture is a single feedforward neural network, for the genus 2 case it is two copies of identical networks, one for each connected component of the disconnected fundamental domain.
The architecture takes the general form
\begin{equation}
    \phi_{\theta} : \Omega \subset \mathbb{R}^2 \to \mathbb{R}^3,
\end{equation}
which maps parameter-space coordinates $(u,v)$ to the embedded point $\phi_{\theta}(u,v) = (x,y,z)$. 

Importantly, a spectral action is used on the fundamental domain inputs to intrinsically enforce the required identification relations; these change dependent on the problem to reflect the surface's topological genus.
The input dimension varies to reflect the basis of these spectral actions, as well as the specified order. 
For genuses $\{0, 1, 2\}$ the input dimensions are $\{(L+1)^2=9, 4N=8, 4N=24\}$ respectively (note genus 2 has two of these networks).
The output is always 3d for the embedding coordinates of the input point in $\mathbb{R}^3$.

The sizes of the dense layers were $(64, 128, 128, 64)$ for genus 0 and 1, then $(128, 256, 512, 256, 128)$ for genus 2.
The layers used tanh activation, which importantly makes the entire architecture a smooth function. 
Genuses 0 and 1 pretrained for 20 epochs with 4096 new points per epoch in batches of size 512, genus 2 did the same for 200 epochs.
Supervised pretraining equally weighted the positional and derivative components, $1 = \beta_1 = \beta_2(t) \ \forall t$.
Then, for the main unsupervised PINN-style training for each genus the training data was resampled every epoch, using (number of epochs, number of points, batch size): genus 0 (30, 1000, 100), genus 1 (200, 4000, 400), genus 2 (2000, 5000, 1000) respectively.
An AdamW optimiser was used with $10^{-5}$ weight decay, and gradient clipping at 1.

The losses used varied scheduling for the total learning rate.
For genus 0 there was none giving a fixed learning rate of 0.003, then for genus 1 and 2 a cosine schedule was used to attune between learning rates (0.00003 to 0.00000001) and (0.0002 to 0.000001) respectively.
The loss components had weights: $\lambda_\mathcal{W}=1, \lambda_R=10, \lambda_G=200$.
Within $\mathcal{L}_R$ the subweightings were: $\alpha_{\text{area}}= 1,\alpha_{\text{pos}}= 0.5,\alpha_{\text{smooth}}= 0.5$, and $\alpha_{\text{log}}= 2$ for the genus 2 case; these used $a_{\text{min}}= 0.01, \epsilon_{\text{pos}}= 0.001, g_{\text{max}}= 5$.
If the annulus regularity weight turns on for genus 2, (radius dropping below 2.5) it uses weight 20 to prevent collapse.
In the genus 2 case there was also a linear warm up schedule for the Willmore loss attuning between weights of 0.01 and 1 over 200 epochs.
Then also for genus 2 the $\alpha_{C^1}$ is zero until 50 epochs where it then ramps linearly over 10 epochs, and equivalently for $\alpha_{C^1}$ ramping over 20 epochs from epoch 150; as these turn on they renormalise the $\lambda_G$ accordingly, which overall decays to half its value over the remaining training time.

In addition to the main training functionality, embedding visualisations and independent testing can be completed with respective visualisations scripts,.
These include functionality for visualising the analytic reference setups used for pretraining of each genus for a variety of parameters, and then also equivalently for the supervised trained embedding models on these analytic surfaces, as used for initialisation.

Furthermore, the codebase also includes functionality for regularity violation rollbacks, such that previously saved epochs can be reinstated if there is a serious regularity error, although this was explicitly not used in any presented runs.
Additionally, there is also functionality to switch the network immersion map to predict the relative position vector to the analytic reference model, instead of the pure direct 3d position of each point, this was again not used and expressed expected limitations with moving surfaces far from the reference surface.

\section{Neural flow results}\label{sec:results}
Details of the training runs for the architectures described in §\ref{sec:architecture} are given below for each considered genus $\{0, 1, 2\}$.
The genus 0 case validates the standard exercise of the round sphere minimising Willmore energy.
The genus 1 case validates the Clifford torus minima, as the global minima for genus 1 surfaces, as proved in the acclaimed work \cite{MarquesNeves2014}.
Finally, the genus 2 case investigates an unproven territory, providing a means to identify minima for more complicated topologies.

\subsection{Genus 0}
The genus 0 functionality allows a choice of any $(a,b,c)$ parameters in defining the initialisation ellipsoid; for the run shown here, values of $(a,b,c) = (2, 1, 0.5)$ were used.
The architecture was first supervised pretrained against this initialisation, then trained in an unsupervised manner with the aim of minimising Willmore energy whilst maintaining regularity.

\begin{figure}[H]
    \centering
    \begin{subfigure}[t]{0.32\textwidth}
        \centering
        \includegraphics[width=\linewidth]{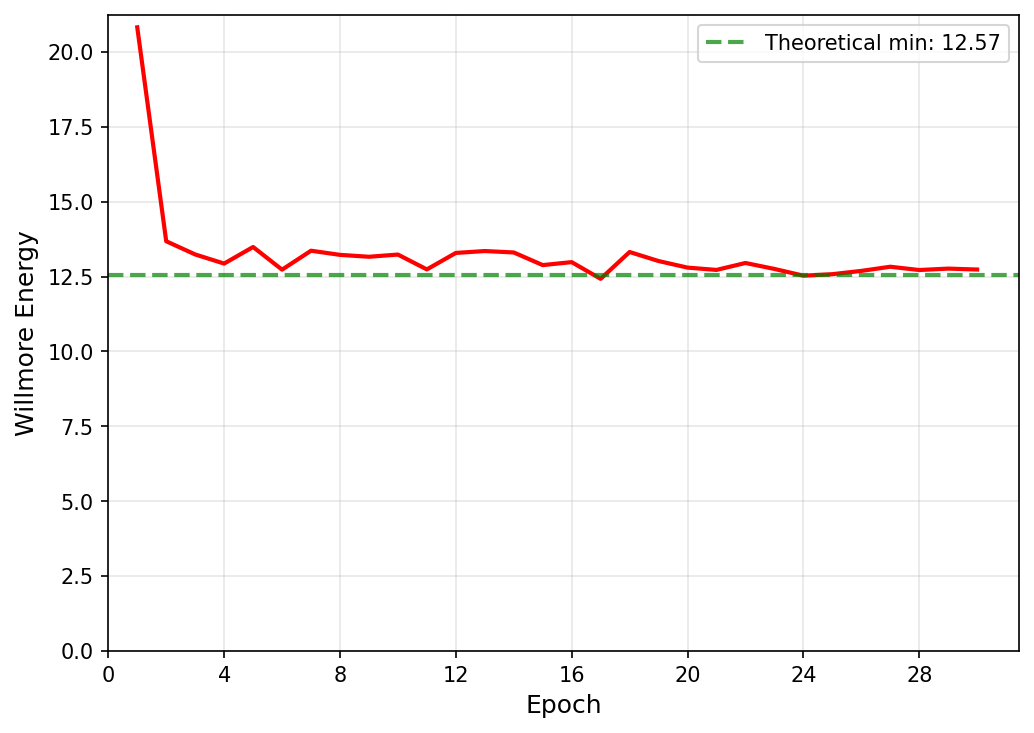}
        \caption{Willmore Loss}
        \label{fig:g0_willmore}
    \end{subfigure}
    \hfill
    \begin{subfigure}[t]{0.32\textwidth}
        \centering
        \includegraphics[width=\linewidth]{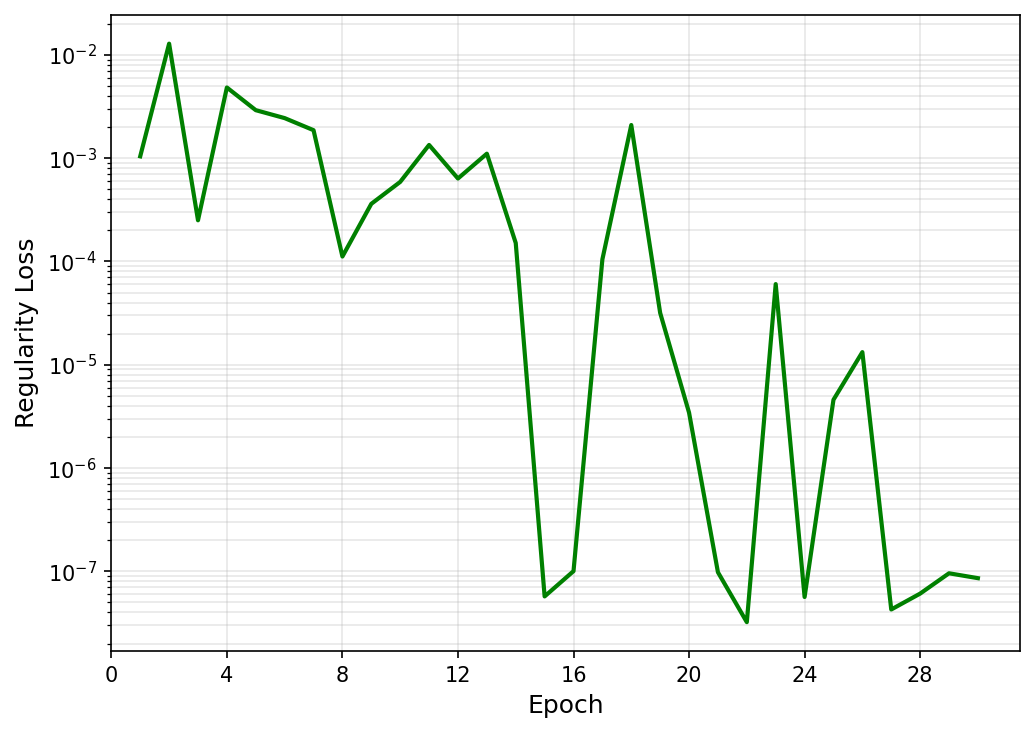}
        \caption{Regularity Loss}
        \label{fig:g0_regularity}
    \end{subfigure}
    \hfill
    \begin{subfigure}[t]{0.32\textwidth}
        \centering
        \includegraphics[width=\linewidth]{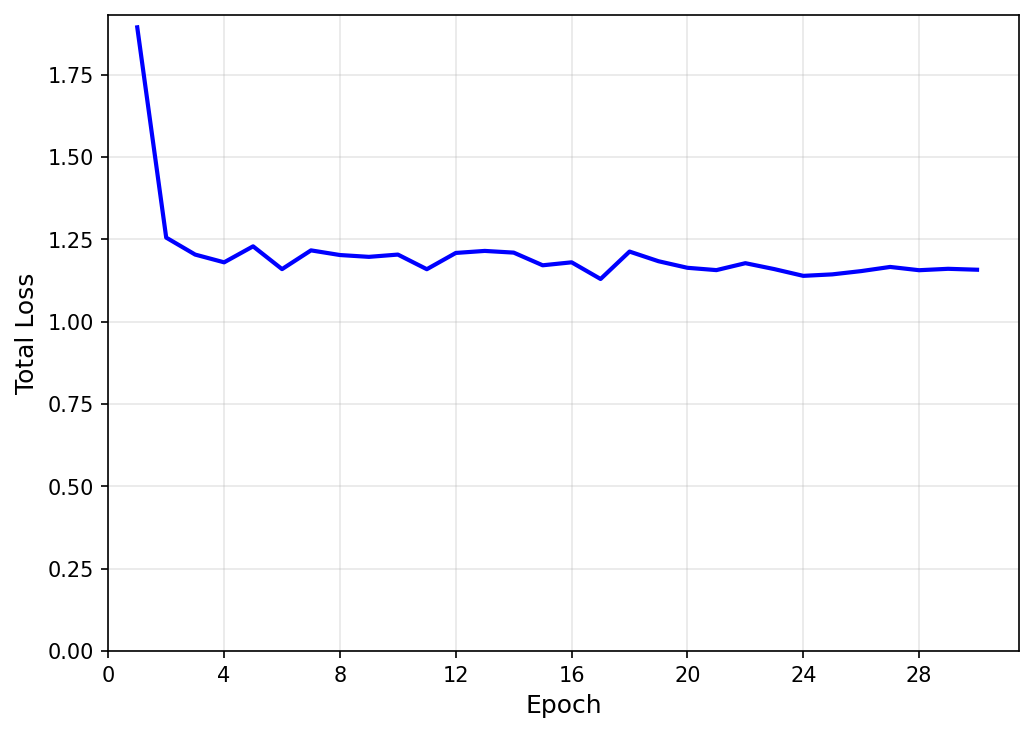}
        \caption{Total Loss}
        \label{fig:g0_total}
    \end{subfigure}
    \caption{Training losses for the genus 0 setup. Willmore loss reaches the expected theoretical minimum for the round sphere; regularity loss stays exceptionally low throughout (far lower order of magnitude) affirming a smooth consistent surface.}
    \label{fig:g0_losses}
\end{figure}

Losses over the 30 epoch training run are shown in Figure \ref{fig:g0_losses}. 
The regularity loss stays negligibly low (note far low orders on the y axis), barely above the ReLU-gating thresholds for each term, this confirms the surface throughout training was smooth and regular. 
Due to these low values the total loss hence reflects the behaviour of solely the Willmore loss.
Then this Willmore loss drops nicely early on, converging towards the theoretical minimum expected from the round sphere.

\begin{figure}[H]
    \centering
    \begin{subfigure}[t]{0.32\textwidth}
        \centering
        \includegraphics[width=0.8\linewidth]{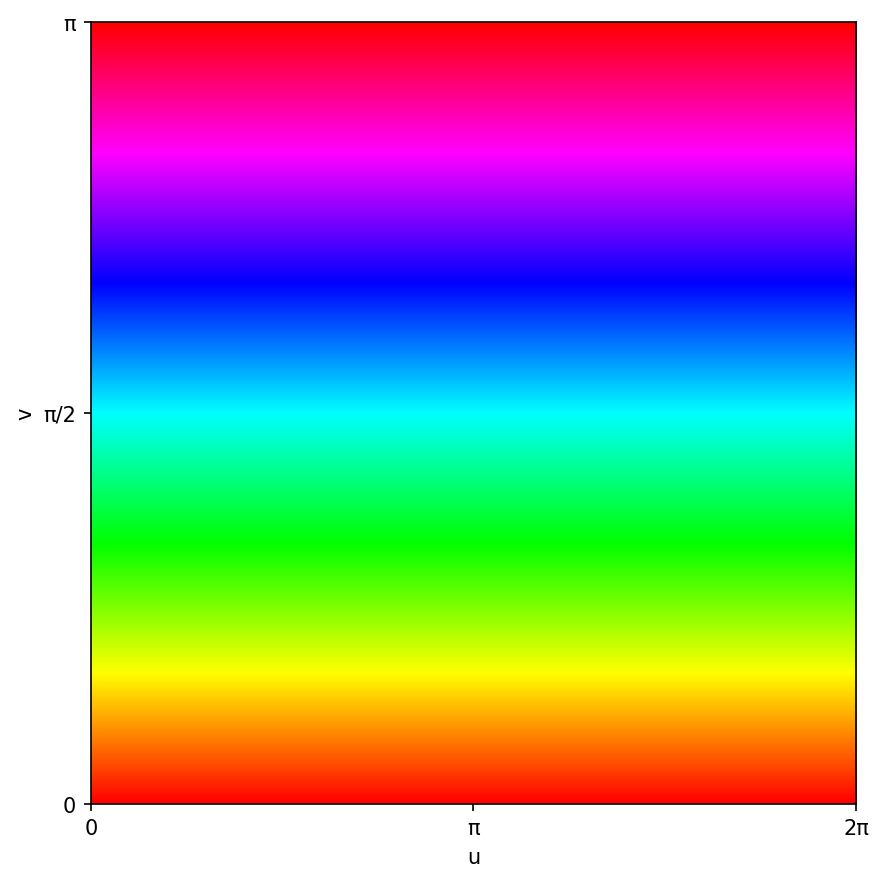}
        \caption{Fundamental Domain}
        \label{fig:g0_fd}
    \end{subfigure}
    \hfill
    \begin{subfigure}[t]{0.32\textwidth}
        \centering
        \includegraphics[width=\linewidth]{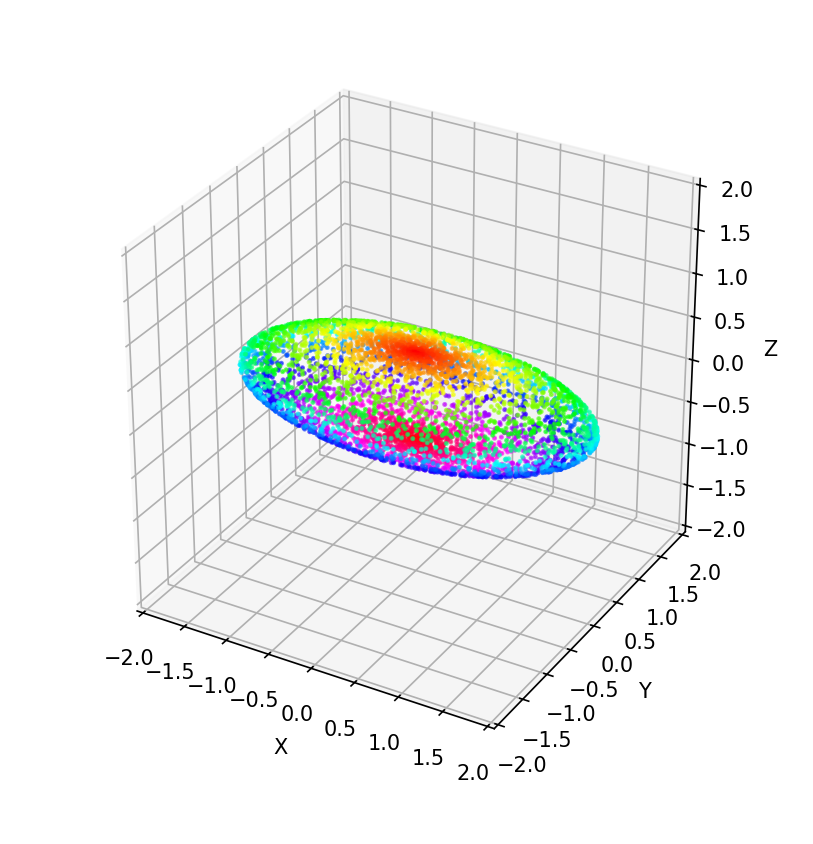}
        \caption{Epoch 0, $\widehat{\mathcal{W}} = 25.62$}
        \label{fig:g0_e0}
    \end{subfigure}
    \hfill
    \begin{subfigure}[t]{0.32\textwidth}
        \centering
        \includegraphics[width=\linewidth]{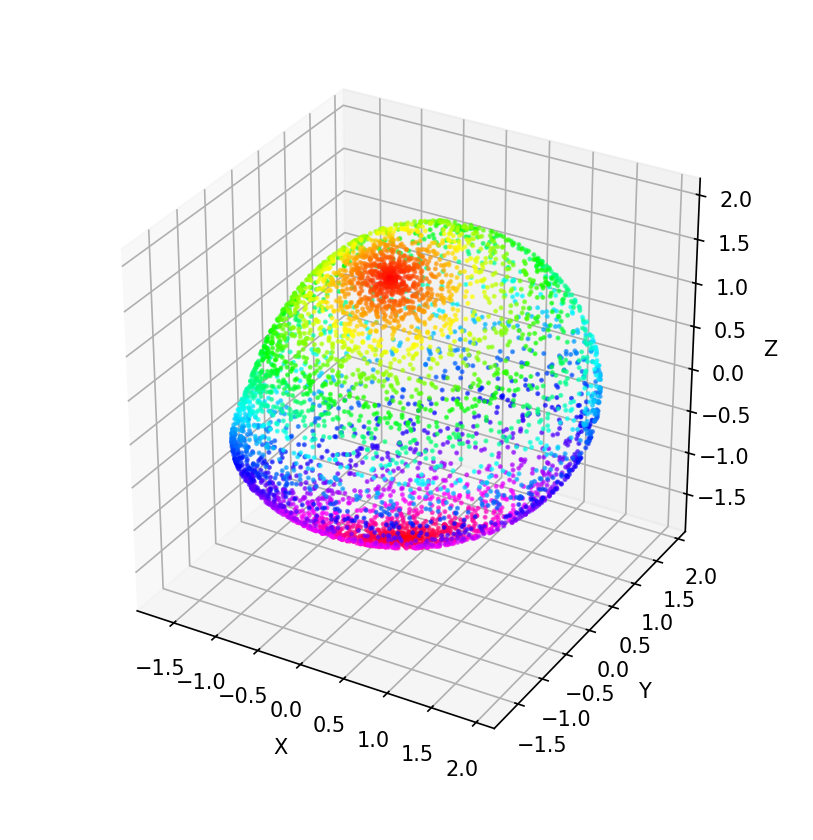}
        \caption{Epoch 1, $\widehat{\mathcal{W}} = 20.82$}
        \label{fig:g0_e1}
    \end{subfigure}\\
    \begin{subfigure}[t]{0.32\textwidth}
        \centering
        \includegraphics[width=\linewidth]{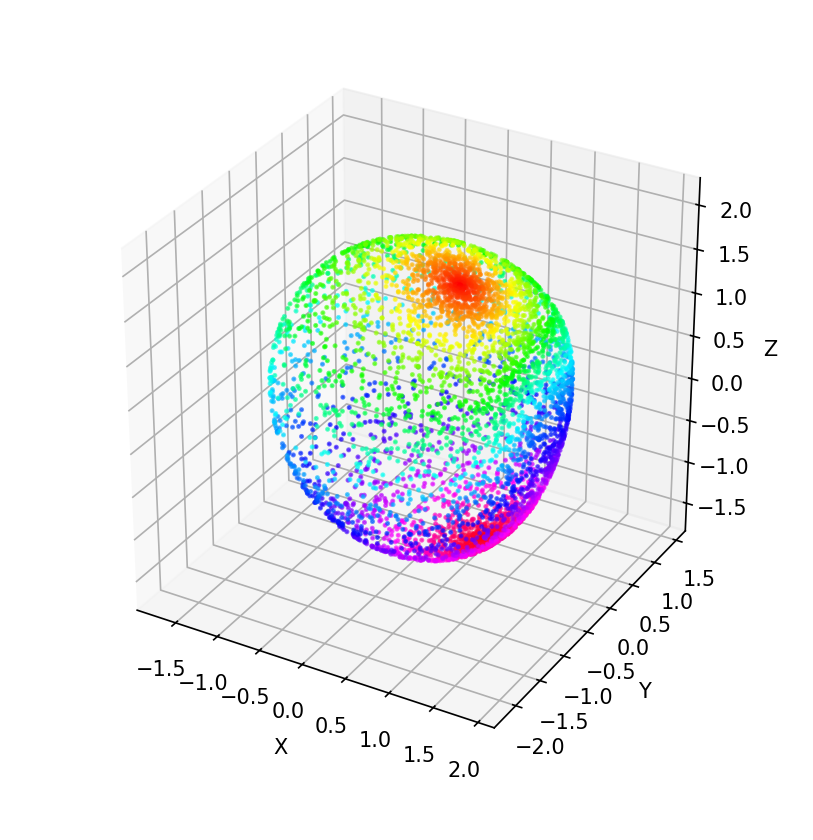}
        \caption{Epoch 6, $\widehat{\mathcal{W}} = 13.48$}
        \label{fig:g0_e6}
    \end{subfigure}
    \hfill
    \begin{subfigure}[t]{0.32\textwidth}
        \centering
        \includegraphics[width=\linewidth]{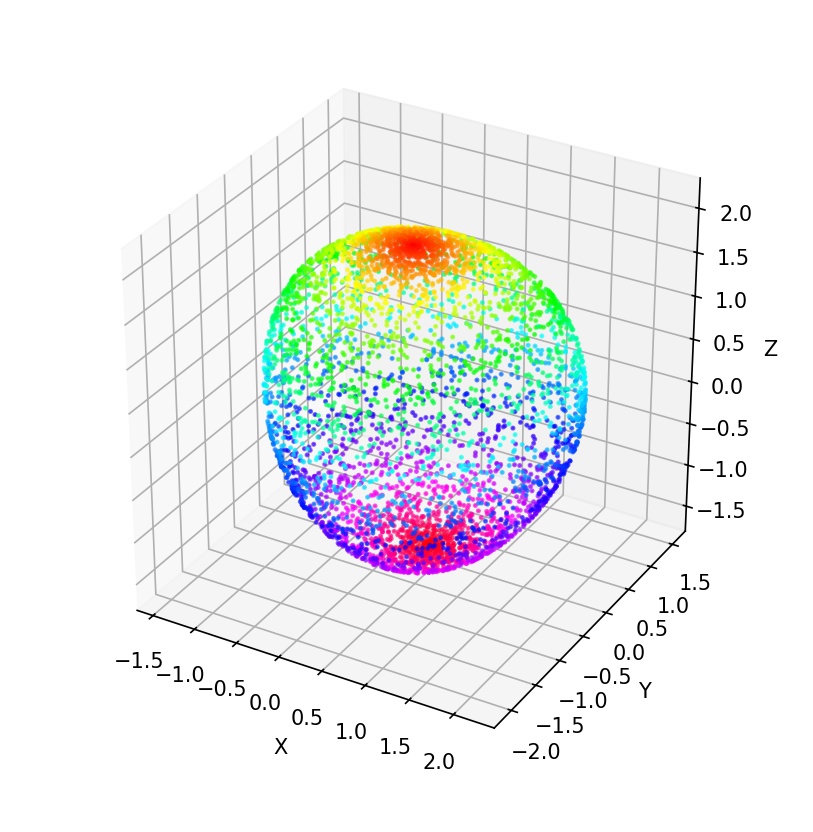}
        \caption{Epoch 15, $\widehat{\mathcal{W}} = 12.88$}
        \label{fig:g0_e15}
    \end{subfigure}
    \hfill
    \begin{subfigure}[t]{0.32\textwidth}
        \centering
        \includegraphics[width=\linewidth]{figures/genus0/embedding_epoch_30.png}
        \caption{Epoch 30, $\widehat{\mathcal{W}} = 12.73$}
        \label{fig:g0_e30}
    \end{subfigure}
    \caption{Visualisation of the embedding evolution over training for the genus 0 setup. The fundamental domain is shown in (a) where the $(u,v)$ coordinates denote the architecture input (after conversion to spherical harmonics), then (b)-(f) are embedding outputs of the model at various epochs through training, using point colours matching the fundamental domain input, with values of the Willmore energy $\widehat{\mathcal{W}}$.}
    \label{fig:g0_visuals}
\end{figure}

Figure \ref{fig:g0_visuals} then shows first a colouring of the fundamental domain of the sphere, which defines the $(u,v)$ coordinates used as input (before passing to the spherical harmonic basis).
A new independent sample of 5000 points in this domain is then taken and passed through the architecture at each stage of training to produce images of the embeddings at each stage.
These embedded points are then plotted in $\mathbb{R}^3$, respecting the fundamental domain colouring, and give the remaining plots in the figure.
These show how the surface is changing over this neural Willmore flow to approach this expected round sphere minimum.
Monte-Carlo estimates of the Willmore energy, $\widehat{\mathcal{W}}$, are also taken on this sample and are reported in the captions, showing a converging decrease towards the expected minimum.

This genus 0 case acts as strong validation for the architecture, and this novel neural approach to modelling the Willmore flow, corroborating the expected convergence to the round sphere surface with Willmore energy $4\pi \sim 12.57$.

\subsection{Genus 1}\label{sec:results_g1}
The genus 1 functionality allows a choice of any complex $\tau$ (with Im$(\tau) > 0$), defining the initialisation torus; for the run shown here $\tau = 0.1i$ was used (no domain shearing / twisting, an example with this is shown in Appendix §\ref{app:g1_twisted}).
The architecture was first supervised pretrained against this initialisation, then again trained in an unsupervised manner with the aim of minimising Willmore energy whilst maintaining regularity.

\begin{figure}[t]
    \centering
    \begin{subfigure}[t]{0.32\textwidth}
        \centering
        \includegraphics[width=\linewidth]{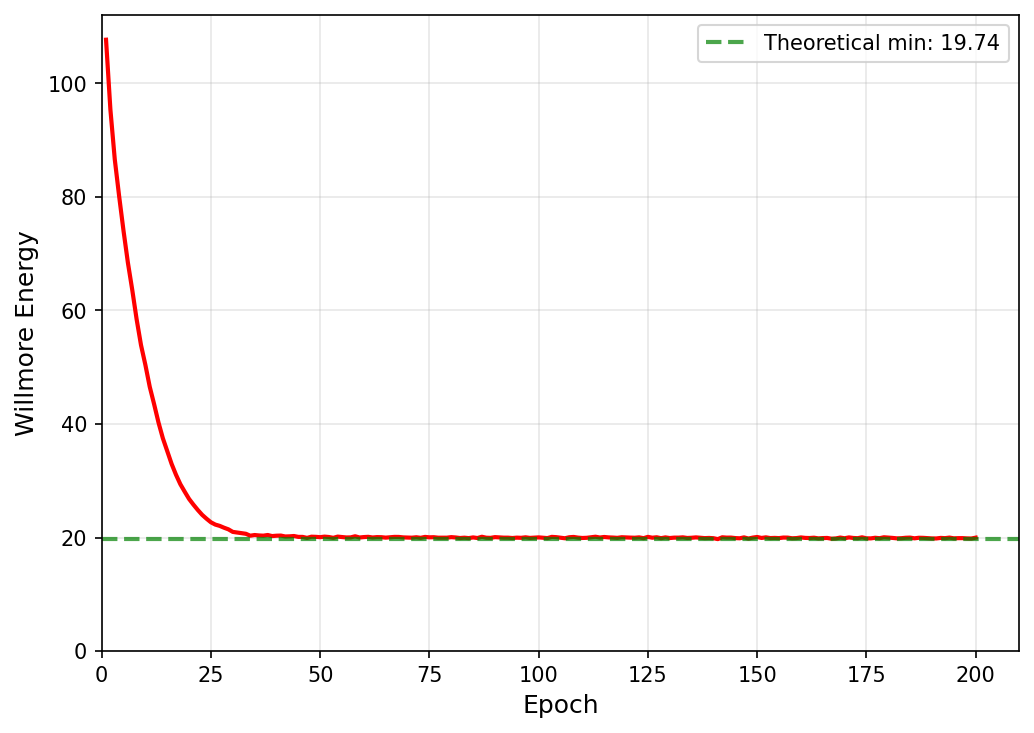}
        \caption{Willmore Loss}
        \label{fig:g1_willmore}
    \end{subfigure}
    \hfill
    \begin{subfigure}[t]{0.32\textwidth}
        \centering
        \includegraphics[width=\linewidth]{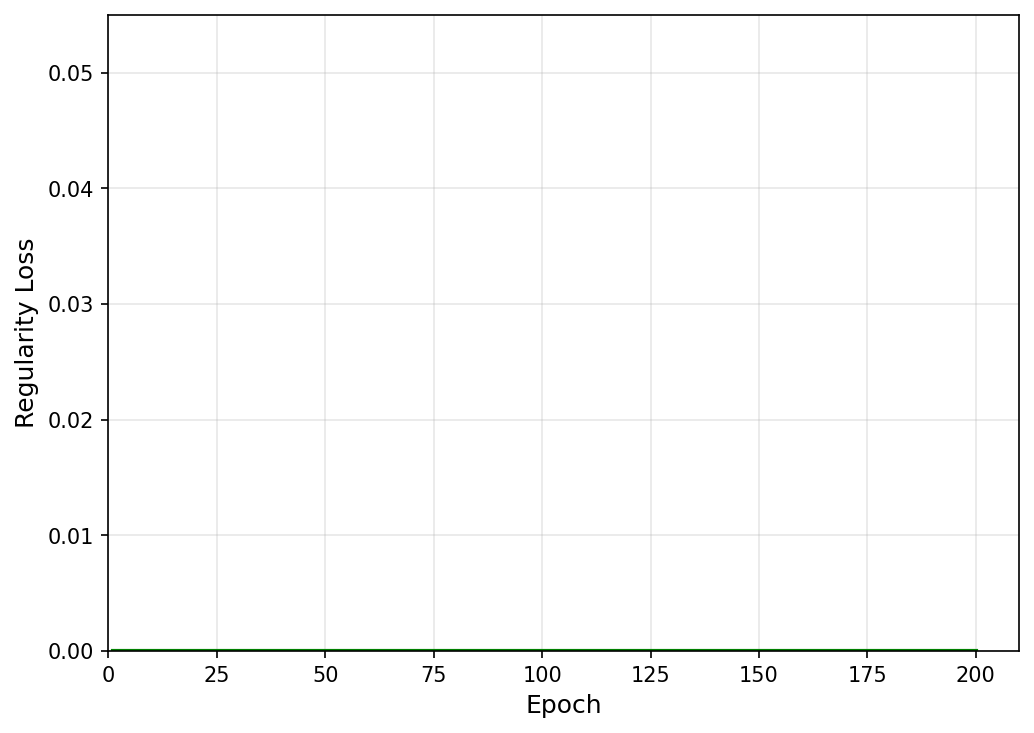}
        \caption{Regularity Loss}
        \label{fig:g1_regularity}
    \end{subfigure}
    \hfill
    \begin{subfigure}[t]{0.32\textwidth}
        \centering
        \includegraphics[width=\linewidth]{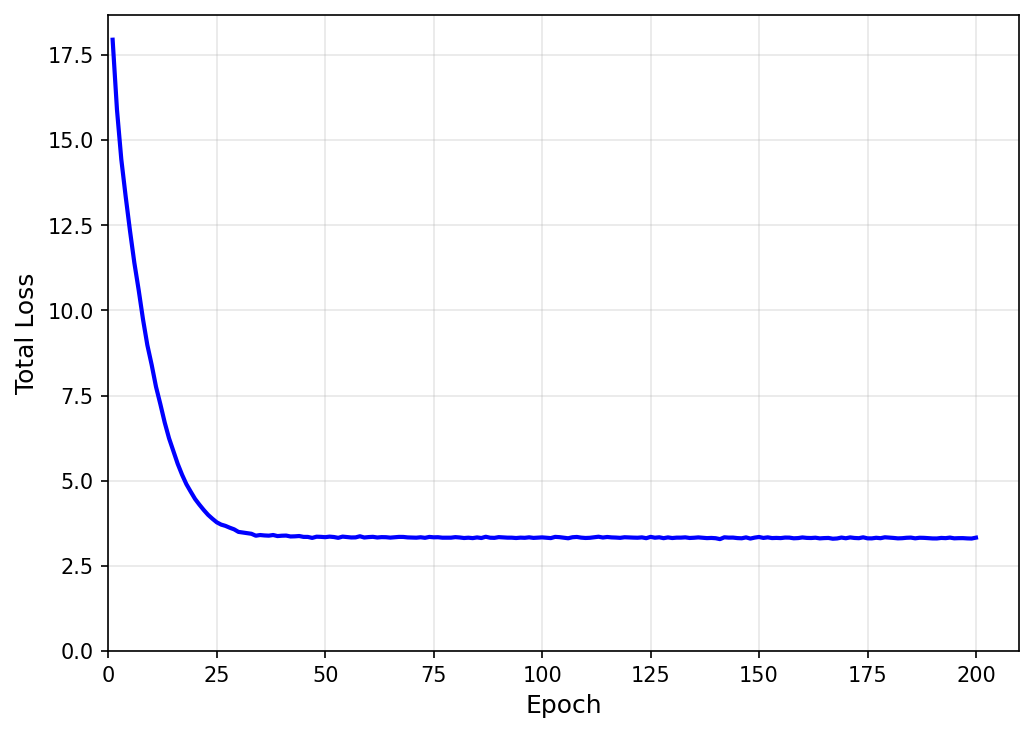}
        \caption{Total Loss}
        \label{fig:g1_total}
    \end{subfigure}
    \caption{Training losses for the genus 1 setup. Willmore loss reaches the expected theoretical minimum for the Clifford torus from the $\tau = 0.1i$ start point; regularity loss stays at exactly 0 throughout (below the ReLU-gated thresholds) affirming a smooth consistent surface.}
    \label{fig:g1_losses}
\end{figure}

Losses over the 200 epoch training run are shown in Figure \ref{fig:g1_losses}. 
In this run the regularity loss stays identically zero, never breaching the ReLU-gating thresholds, asserting that the surface throughout training was smooth and regular. 
The total loss shown hence reflects exactly the behaviour of the Willmore loss.
This Willmore loss drops nicely\footnote{The curve is smoother than genus 0 since more data is used per epoch, and training runs for more epochs.} early on, converging towards this theoretical minimum expected for the Clifford torus.

\begin{figure}[H]
    \centering
    \begin{subfigure}[t]{0.32\textwidth}
        \centering
        \includegraphics[width=0.8\linewidth]{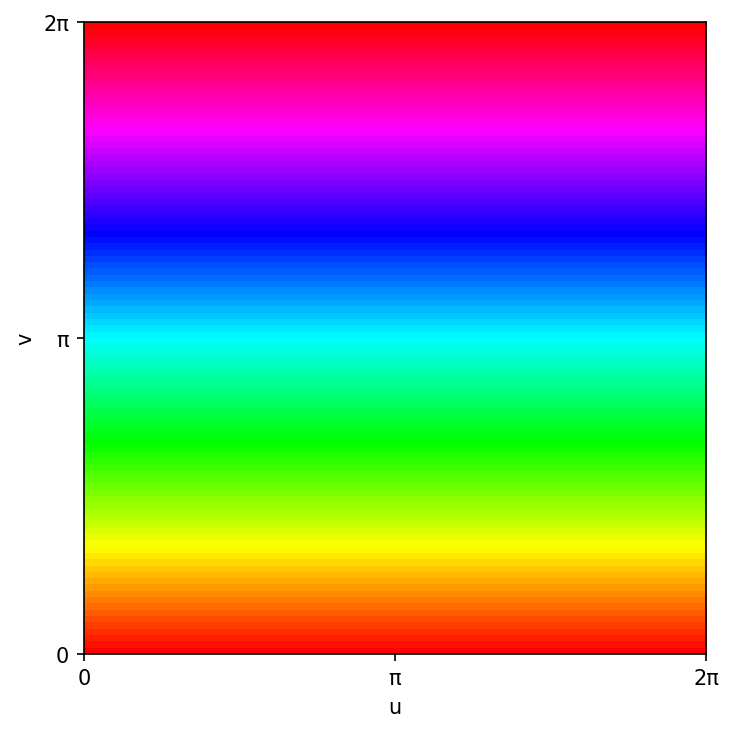}
        \caption{Fundamental Domain}
        \label{fig:g1_fd}
    \end{subfigure}
    \hfill
    \begin{subfigure}[t]{0.32\textwidth}
        \centering
        \includegraphics[width=\linewidth]{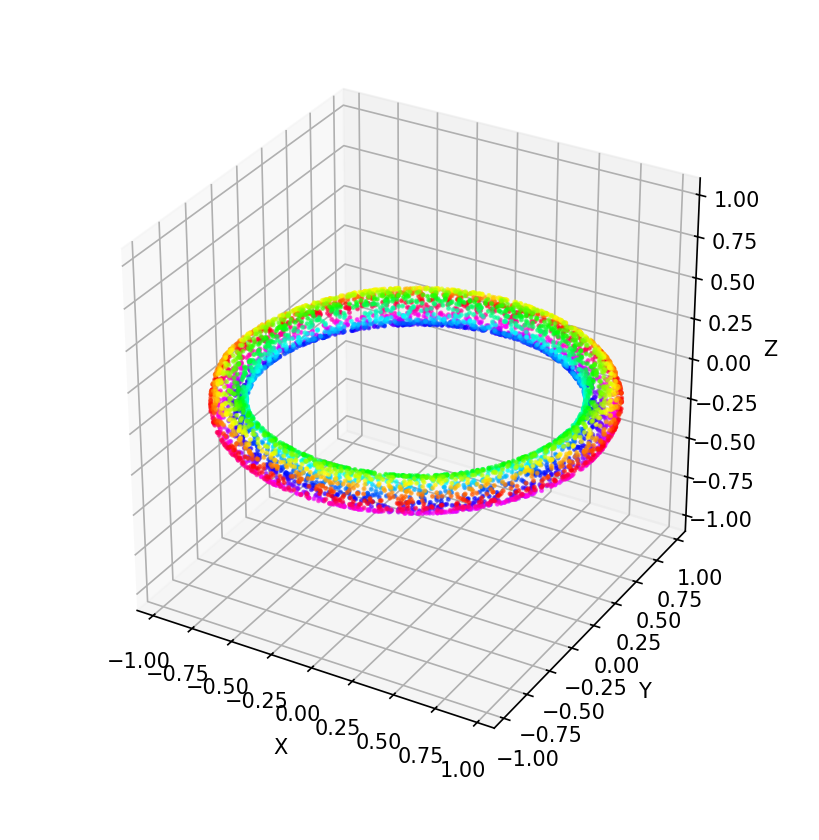}
        \caption{Epoch 0, $\widehat{\mathcal{W}} = 118.97$}
        \label{fig:g1_e0}
    \end{subfigure}
    \hfill
    \begin{subfigure}[t]{0.32\textwidth}
        \centering
        \includegraphics[width=\linewidth]{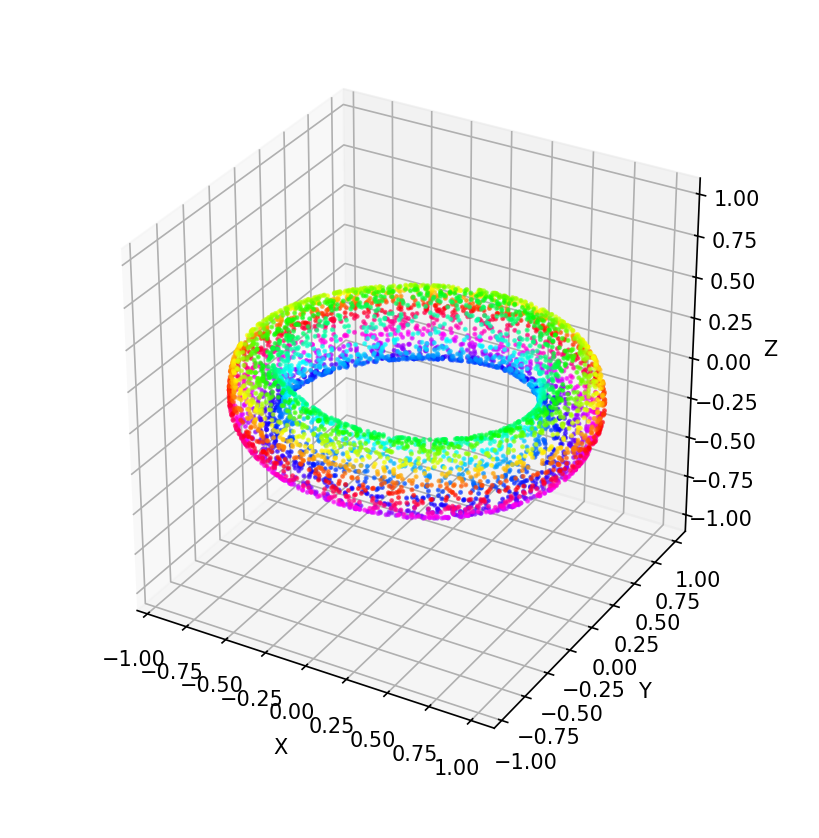}
        \caption{Epoch 10, $\widehat{\mathcal{W}} = 50.42$}
        \label{fig:g1_e10}
    \end{subfigure}\\
    \begin{subfigure}[t]{0.32\textwidth}
        \centering
        \includegraphics[width=\linewidth]{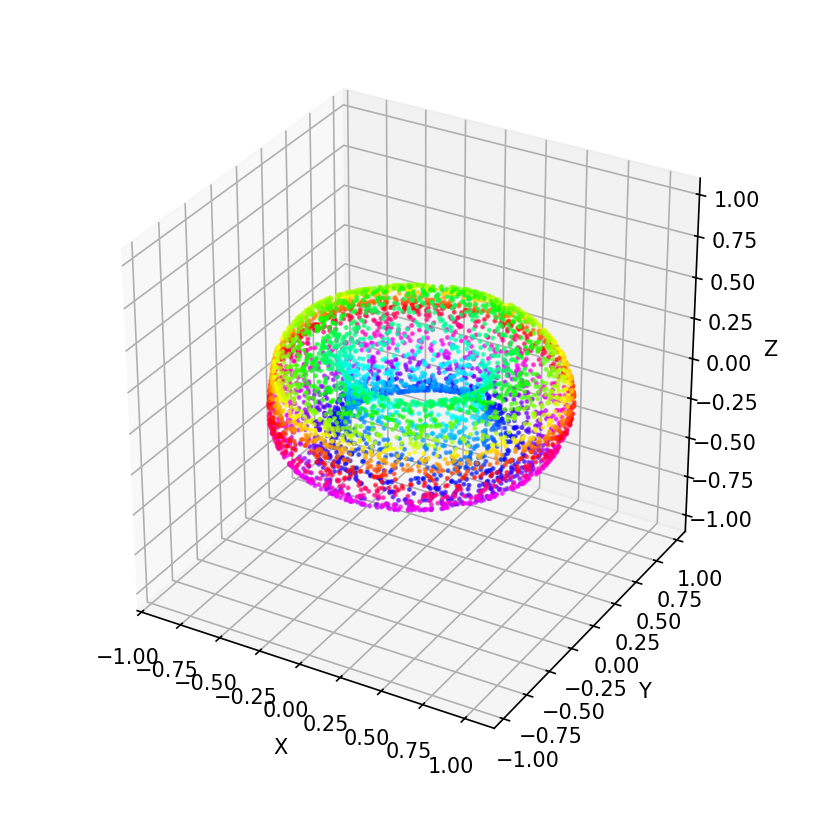}
        \caption{Epoch 20, $\widehat{\mathcal{W}} = 26.80$}
        \label{fig:g1_e20}
    \end{subfigure}
    \hfill
    \begin{subfigure}[t]{0.32\textwidth}
        \centering
        \includegraphics[width=\linewidth]{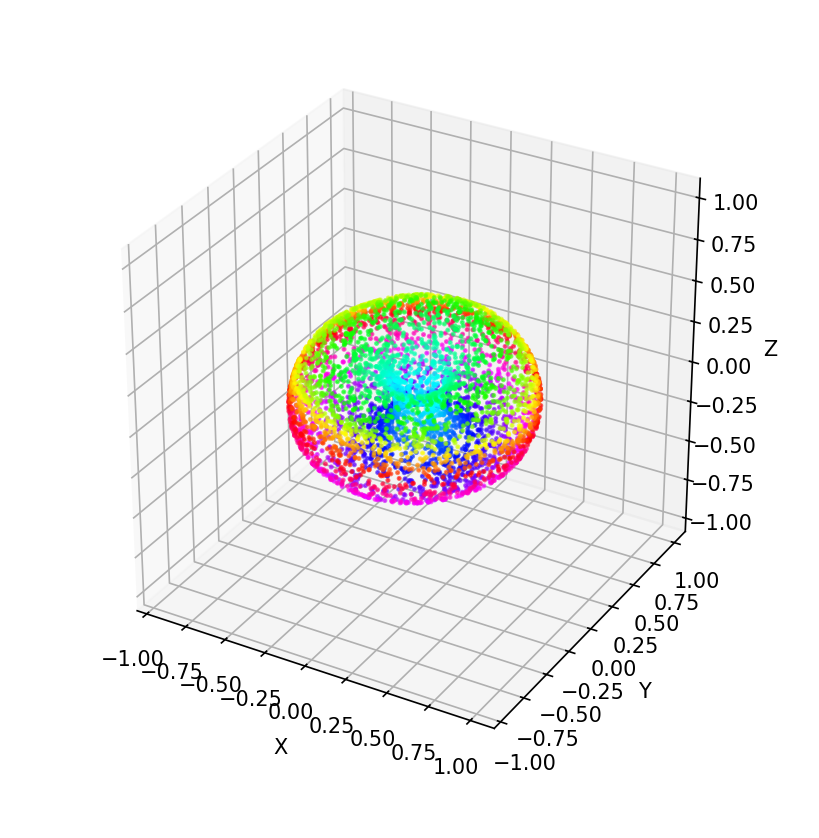}
        \caption{Epoch 100, $\widehat{\mathcal{W}} = 20.03$}
        \label{fig:g1_e100}
    \end{subfigure}
    \hfill
    \begin{subfigure}[t]{0.32\textwidth}
        \centering
        \includegraphics[width=\linewidth]{figures/genus1/embedding_epoch_200.png}
        \caption{Epoch 200, $\widehat{\mathcal{W}} = 19.98$}
        \label{fig:g1_e200}
    \end{subfigure}
    \caption{Visualisation of the embedding evolution over training for the genus 1 setup. The fundamental domain is shown in (a) where the $(u,v)$ coordinates denote the architecture input (after conversion to Fourier modes), then (b)-(f) are embedding outputs of the model at various epochs through training from $\tau = 0.1i$, using point colours matching the fundamental domain input, with values of the Willmore energy $\widehat{\mathcal{W}}$.}
    \label{fig:g1_visuals}
\end{figure}

Figure \ref{fig:g1_visuals} then shows first a colouring of the fundamental domain of the torus, which defines the $(u,v)$ coordinates used as input (before passing to the Fourier basis).
Again, a new sample of 5000 points in this domain is then taken and passed through the architecture at each stage of training to produce images of the embeddings at each stage.
These embedded points are then plotted in $\mathbb{R}^3$, respecting the fundamental domain colouring, and give the remaining plots in the figure.
These embedding plots show how the surface changes over this neural Willmore flow, approaching this expected Clifford torus minimum.
Monte-Carlo estimates of the Willmore energy, $\widehat{\mathcal{W}}$, are also taken on this sample at each stage and are reported in the captions, again showing a converging decrease towards the expected minimum.

The genus 1 case is a truly non-trivial problem, proving it \cite{MarquesNeves2014} required extremely advanced techniques, and is regarded one of the most significant proofs of this century.
Here this novel neural approach to modelling the Willmore flow strongly corroborates this result.
The architecture has complete freedom to construct any immersion map it wants, with no symmetry or intersection constraints, yet when guided by the losses built in this work converges towards this symmetric Clifford torus embedding solution.
Overall, this demonstrates the method can perform well in complicated scenarios such as this convergence to the Clifford torus surface with minimal Willmore energy $2\pi^2 \sim 19.74$.

\subsection{Genus 2}
The genus 2 case is an open problem, with many excellent works examining it from an analytic viewpoint since its conception.
Here, this neural Willmore flow machine learning architecture, previously validated on the genus 0 and genus 1 cases, is applied to this open genus 2 scenario.

The construction for genus 2 surfaces in the codebase is more subtle, using the gluing of two punctured tori.
The genus 2 code functionality allows choice of any complex $\tau_1, \tau_2$ for these two tori domains, and any choice of the puncture hole radius $\delta$, to define this initialisation setup.
Since this is open territory, a more conservative start point is used in the run shown here, with $\tau_1=\tau_2=0.7$, and $\delta=0.65$.
The two tori are hence initialised close to Clifford tori, such that their direct sum of Willmore energies is already low, despite being undefined at the lack of initial intersection.
The architecture was first supervised pretrained against this initialisation, noting that it produces something disconnected in the embedding $\mathbb{R}^3$ space.
The following unsupervised training then first prioritises the gluing loss (sequentially $C^0$, then $C^1$, then $C^2$) to form a consistent smooth surface, before the Willmore loss takes over as the training process looks to minimise the Willmore energy whilst maintaining regularity.

\begin{figure}[b]
    \centering
    \begin{subfigure}[t]{0.32\textwidth}
        \centering
        \includegraphics[width=\linewidth]{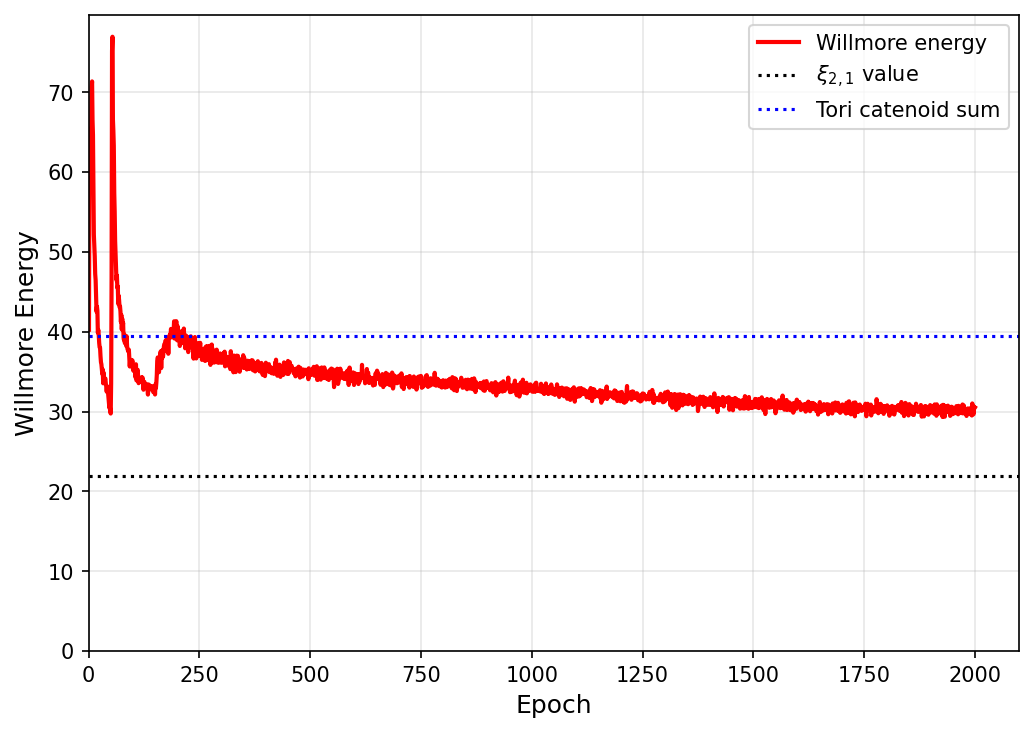}
        \caption{Willmore Loss}
        \label{fig:g2_willmore}
    \end{subfigure}
    \hfill
    \begin{subfigure}[t]{0.32\textwidth}
        \centering
        \includegraphics[width=\linewidth]{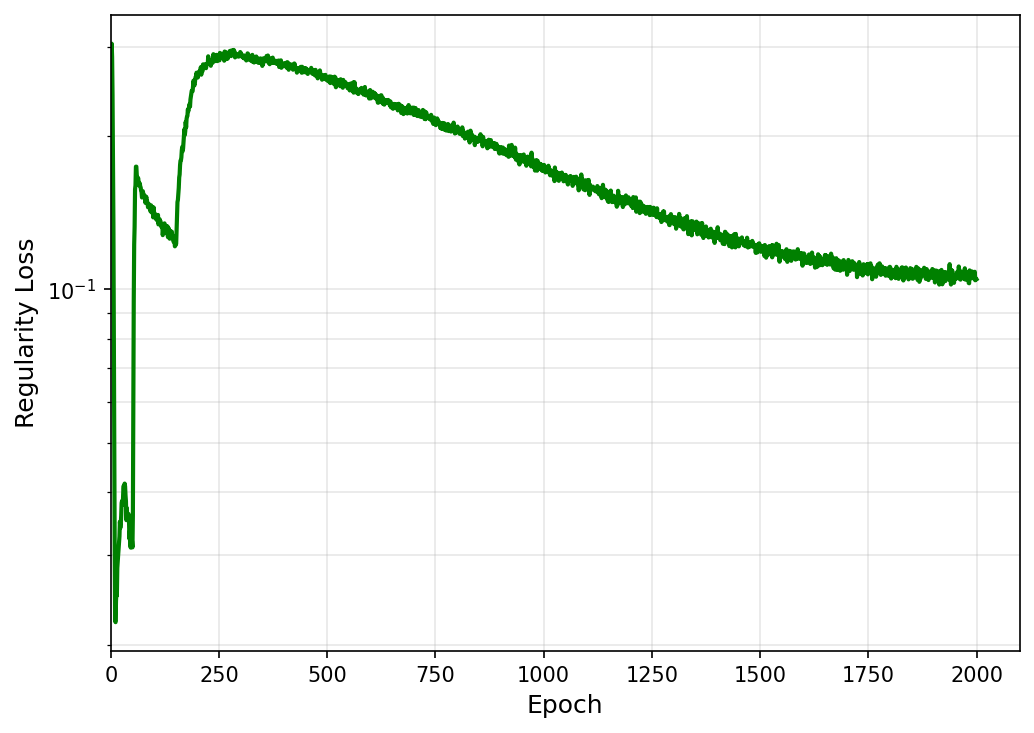}
        \caption{Regularity Loss}
        \label{fig:g2_regularity}
    \end{subfigure}
    \hfill
    \begin{subfigure}[t]{0.32\textwidth}
        \centering
        \includegraphics[width=\linewidth]{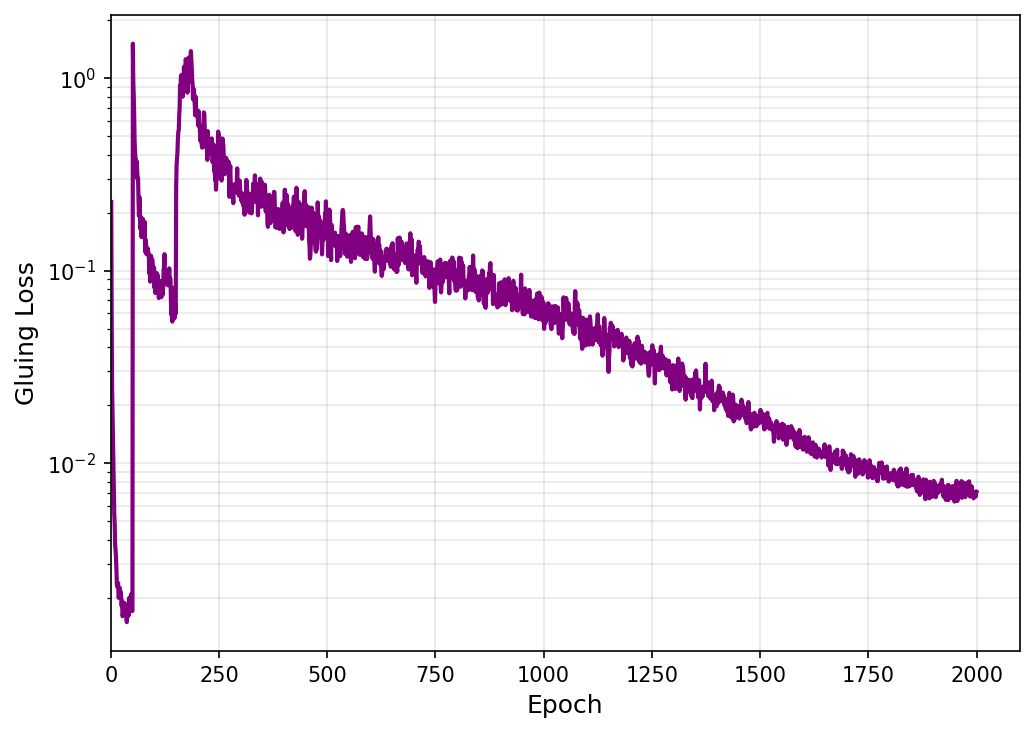}
        \caption{Gluing Loss}
        \label{fig:g2_gluing}
    \end{subfigure}
    \caption{Training losses for the genus 2 setup. Willmore loss has a more erratic decline but drops to a stable region. Regularity loss stays low throughout affirming a smooth consistent surface. Since this setup used a disconnected fundamental domain a gluing loss also was used which was a weighted sum of $(C^0, C^1, C^2)$ components, this also dropped substantially over the training as the gluing junction was smoothened out; at epoch 50 when the $C^1$ component was introduced there is a spike in both gluing and Willmore losses, and the same at epoch 150 when $C^2$ component was introduced as the surface restabilises. Total loss is not shown here since component weight scheduling makes it non-monotonic and less informative.}
    \label{fig:g2_losses}
\end{figure}

Losses over the 2000 epoch training run are shown in Figure \ref{fig:g2_losses}. 
In all the losses the impact of the $C^1$ and then $C^2$ ramp delays, at epochs 50 and 150 respectively, are clearly visible.
Hence behaviour after these peaks is most important.
Due to the varied loss component scheduling the total loss is not monotonic and hence less informative, therefore this figure focuses on each of the loss components directly.

In this run the regularity loss decreases smoothly after the full gluing has been enforced and the surface is consistent.
The gluing loss also shows smooth decrease as the consistency between the tori fundamental domain charts improves throughout the training.
Finally, the Willmore loss also drops nicely, despite the scale being thrown off by the spikes as gluing continuity components are introduced.
This Willmore loss drops throughout this stable part of training, well below the value for the catenoid gluing of two Clifford tori.

\begin{figure}[H]
    \centering
    \begin{subfigure}[t]{0.32\textwidth}
        \centering
        \includegraphics[width=0.95\linewidth]{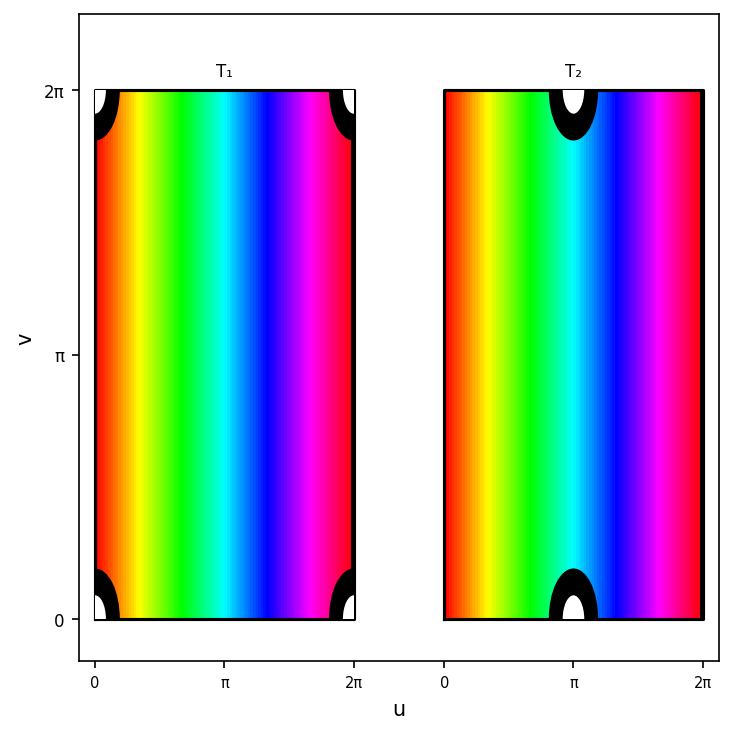}
        \caption{Fundamental Domain}
        \label{fig:g2_fd}
    \end{subfigure}
    \hfill
    \begin{subfigure}[t]{0.32\textwidth}
        \centering
        \includegraphics[width=\linewidth]{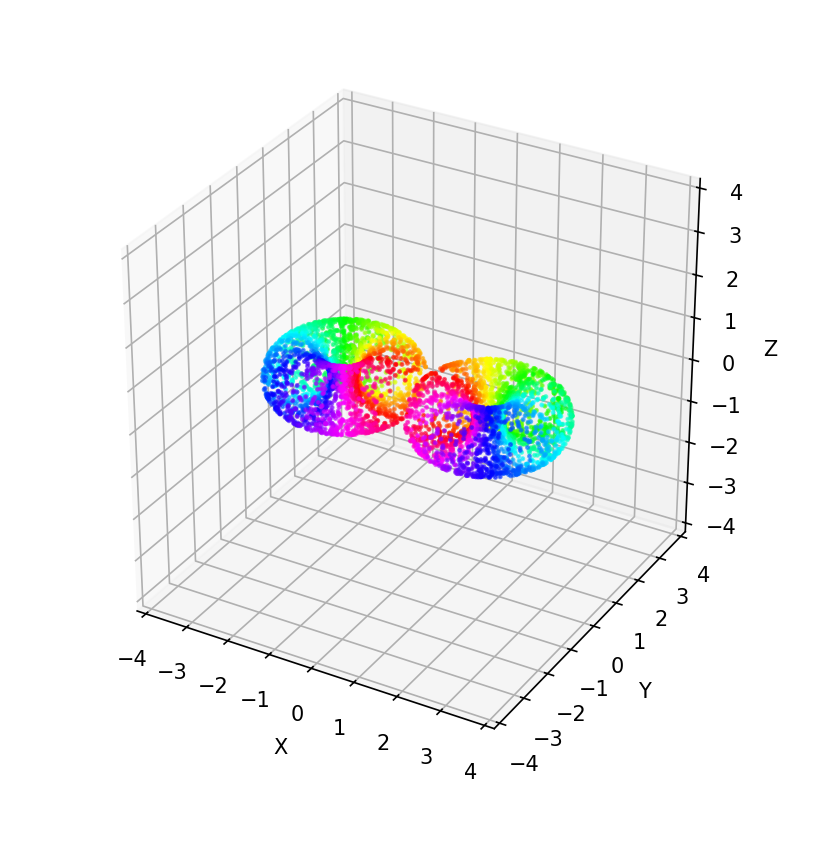}
        \caption{Epoch 0, $\widehat{\mathcal{W}}^* = 36.59$}
        \label{fig:g2_e0}
    \end{subfigure}
    \hfill
    \begin{subfigure}[t]{0.32\textwidth}
        \centering
        \includegraphics[width=\linewidth]{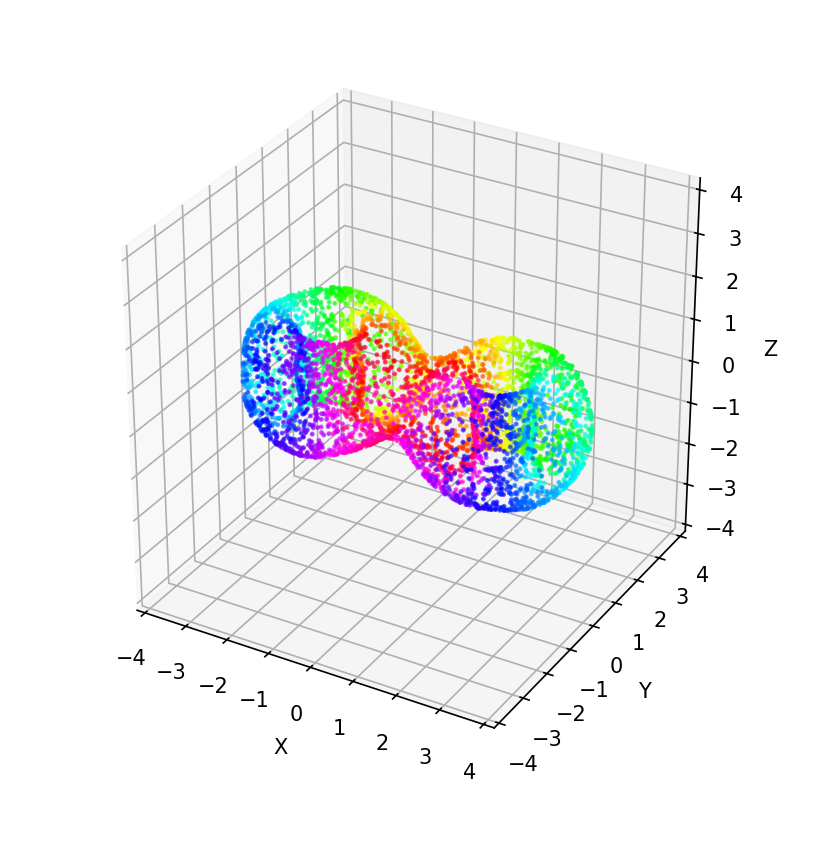}
        \caption{Epoch 10, $\widehat{\mathcal{W}}^* = 50.48$}
        \label{fig:g2_e10}
    \end{subfigure}\\
    \begin{subfigure}[t]{0.32\textwidth}
        \centering
        \includegraphics[width=\linewidth]{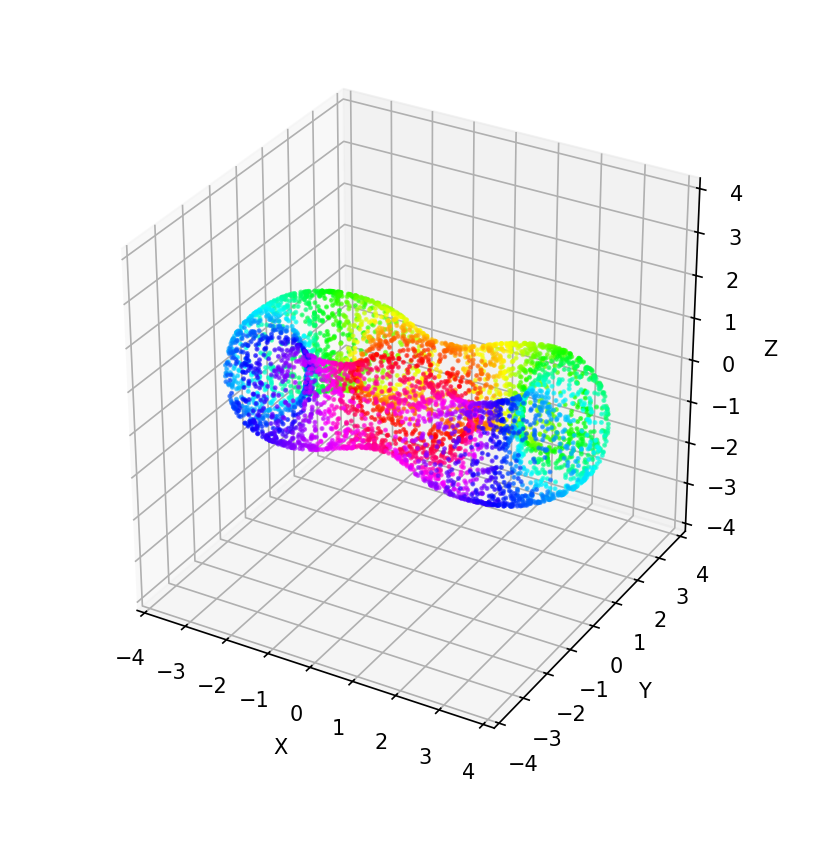}
        \caption{Epoch 25, $\widehat{\mathcal{W}}^* = 34.73$}
        \label{fig:g2_e25}
    \end{subfigure}
    \hfill
    \begin{subfigure}[t]{0.32\textwidth}
        \centering
        \includegraphics[width=\linewidth]{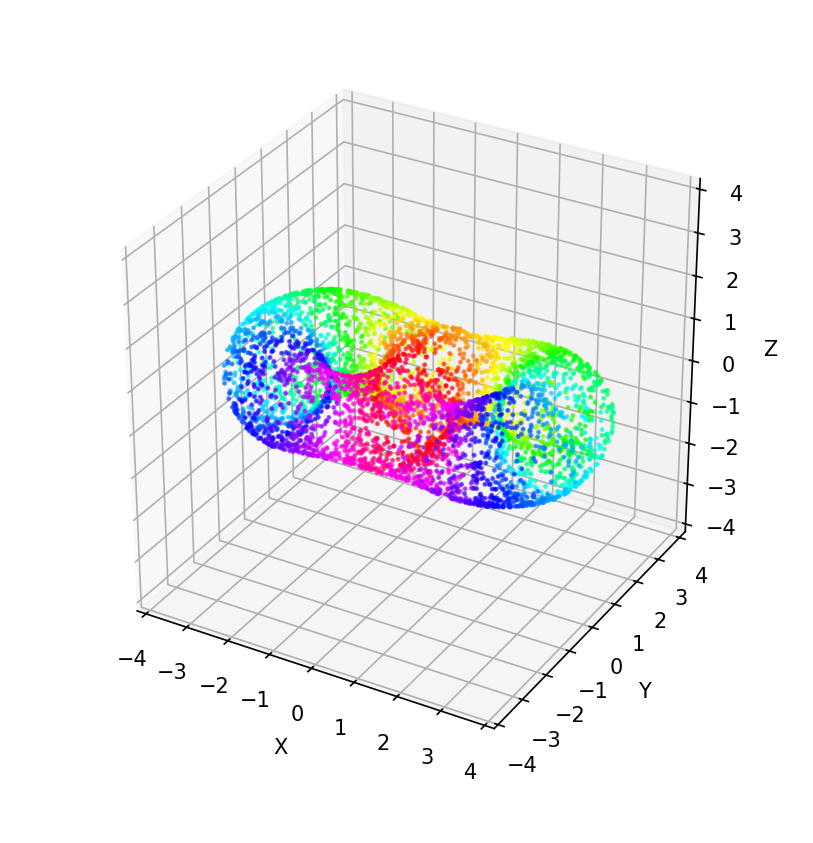}
        \caption{Epoch 50, $\widehat{\mathcal{W}}^* = 28.85$}
        \label{fig:g2_e50}
    \end{subfigure}
    \hfill
    \begin{subfigure}[t]{0.32\textwidth}
        \centering
        \includegraphics[width=\linewidth]{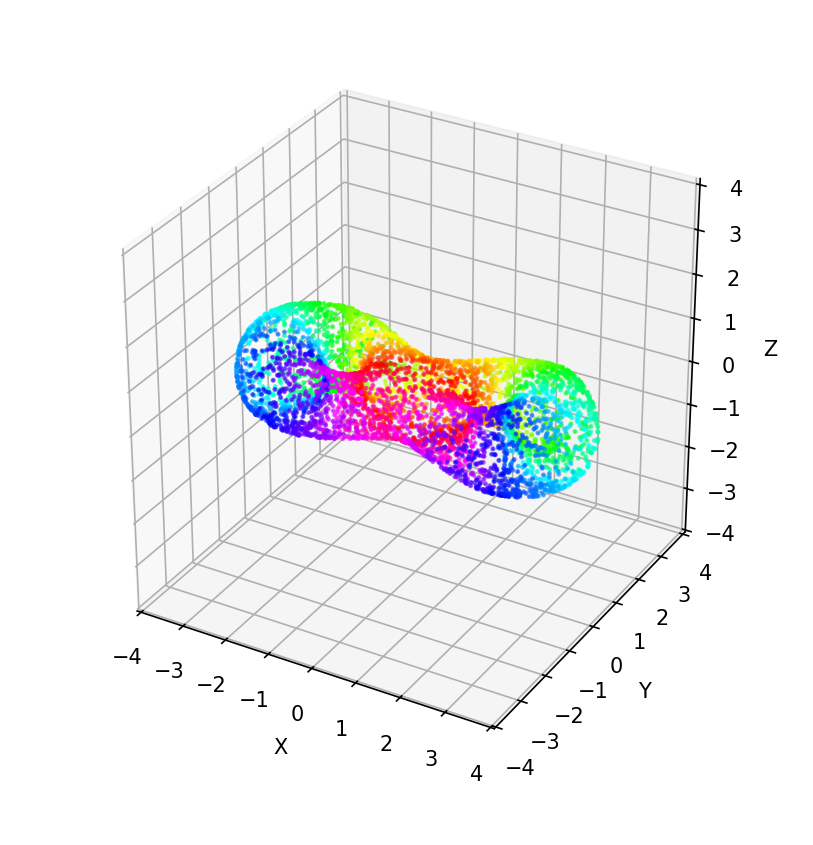}
        \caption{Epoch 100, $\widehat{\mathcal{W}}^* = 35.39$}
        \label{fig:g2_e100}
    \end{subfigure}\\
    \begin{subfigure}[t]{0.32\textwidth}
        \centering
        \includegraphics[width=\linewidth]{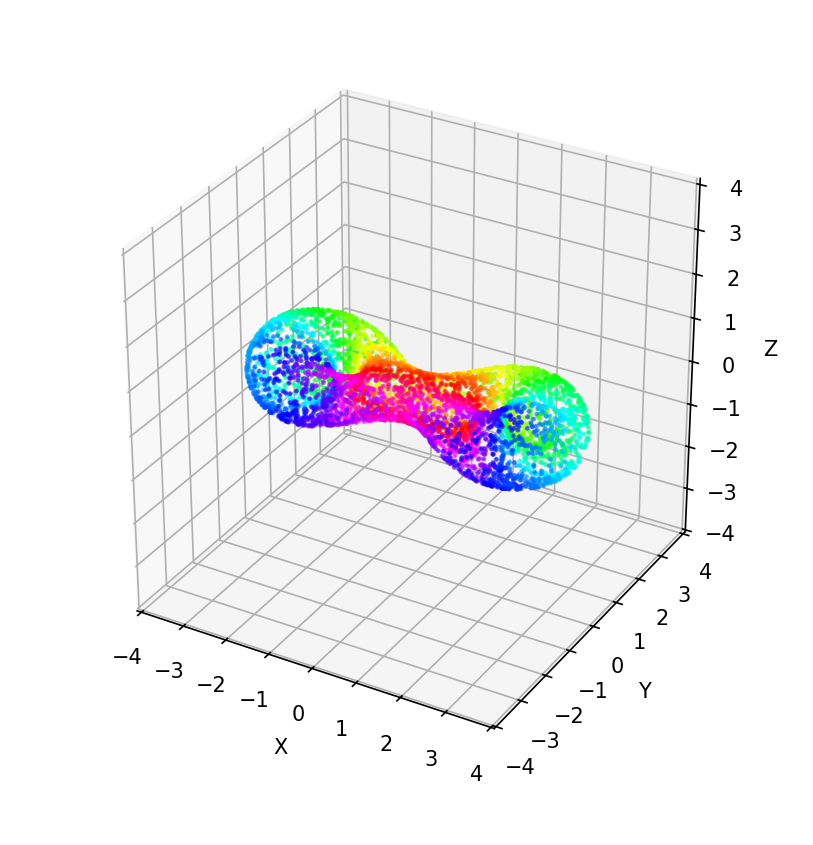}
        \caption{Epoch 500, $\widehat{\mathcal{W}} = 34.82$}
        \label{fig:g2_e500}
    \end{subfigure}
    \hfill
    \begin{subfigure}[t]{0.32\textwidth}
        \centering
        \includegraphics[width=\linewidth]{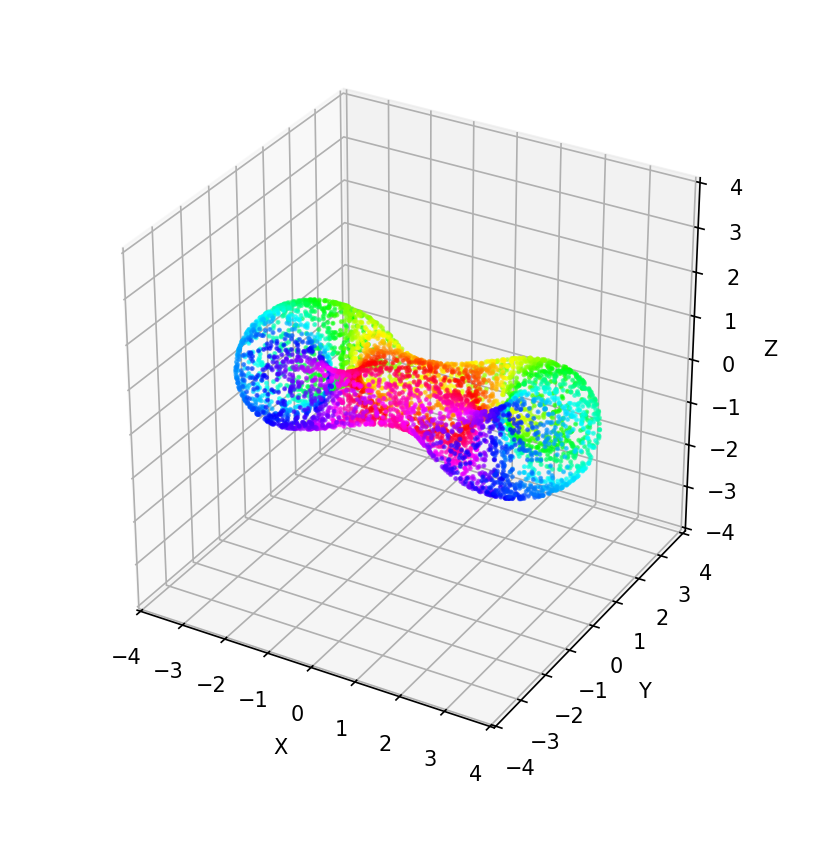}
        \caption{Epoch 1000, $\widehat{\mathcal{W}} = 32.57$}
        \label{fig:g2_e1000}
    \end{subfigure}
    \hfill
    \begin{subfigure}[t]{0.32\textwidth}
        \centering
        \includegraphics[width=\linewidth]{figures/genus2/embedding_epoch_2000.png}
        \caption{Epoch 2000, $\widehat{\mathcal{W}} = 30.19$}
        \label{fig:g2_e2000}
    \end{subfigure}
    \caption{Visualisation of the embedding evolution over training for the genus 2 setup. The fundamentaxl domain is shown in (a) where the $(u,v)$ coordinates denote the architecture input (after conversion to Fourier modes), there're 2 domains for each tori used in gluing. Then (b)-(i) are embedding outputs of the model at various epochs through training, using point colours matching the fundamental domain input. values through training are given for the Willmore energy $\widehat{\mathcal{W}}$, we note a caveat that before $C^2$ gluing is enforced at epoch 150 the measure is unphysical.}
    \label{fig:g2_visuals}
\end{figure}

Figure \ref{fig:g2_visuals} then shows first a colouring of the double fundamental domain set up, which defines the respective $(u,v)$ coordinates used as input for each subnetwork (before passing to the Fourier bases).
Equivalently, a new sample of 5000 points across this domain is taken and passed through the architecture at each stage of training to produce images of the embeddings at each stage.
These embedded points are then plotted in $\mathbb{R}^3$, respecting the fundamental domain colouring, and give the remaining plots in this figure.
These embedding plots show how the surface changes over this neural Willmore flow, where now the initial disconnected behaviour is visible at genus 0, followed by rugged connection which then smoothens out as the $\mathcal{L}_{C^1}$, then $\mathcal{L}_{C^2}$, losses turn on.
Finally, the full surface adjusts as the tori holes move toward the junction and the constituent tori become asymmetric around their domains, producing a smooth surface with low Willmore energy; as reported using Monte-Carlo estimates over the sample.

The genus 2 case is an open area for investigation, and the numerics of this genus-0 and genus-1 validated novel neural approach do show smoothening of a consistent genus 2 surface.
The resulting surface has a numerical Willmore energy of 30.19, notably below the value for the catenoid sum of two Clifford tori, $\sim 39.48$.

The visualisations of the learnt surface are unlikely to represent\footnote{The authors refer to \cite{lawson_visuals} for some interesting visualisations of Lawson surfaces, although with focus on them as embeddings in $S^3$ instead of $\mathbb{R}^3$.} the Lawson surface $\xi_{2,1}$.
The surface starts from two disconnected punctured tori domains, where learning would then connect these punctured boundaries and smooth the junction between them.
Naively one would expect production of a catenoid bridge which would add zero Willmore energy to the sum of the nearly Clifford tori used for initialisation.
However, the learning shows that after connection the tori are also adjusted and the bridge does not look catenoid in nature, as energy is taken out of the tori handles into the bridge in a way that minimises it for the full surface; exactly as desired for the neural Willmore flow learning.

The learnt surface has not yet reached the Willmore energy order of the Lawson surface $\xi_{2,1}$, yet perhaps with more compute this would be achievable. 
It is worth noting also that there are error contributions to the estimated Willmore energies associated to the Monte Carlo numerical estimate, as well as the non-zero (yet still small) gluing and regularity losses, and it would be interesting future work to extend training and refine tuning to see how the Willmore estimate changes with lower values of these losses, and if convergence towards $\xi_{2,1}$ is achievable.
Overall, the production of a consistent genus 2 surface with Willmore energy well below the naive catenoid bridge sum of tori demonstrates how this new AI-inspired numerical approach can learn near-minimal surfaces, and potentially suggest new directions for analytic investigation.

\section{Summary}
This work introduces a novel AI-inspired approach to the search for minimal surfaces, specifically to surfaces of minimal Willmore energy.

An embedding map from a 2d fundamental domain is represented with a neural architecture, taking 2d inputs for points on the domain and mapping them to 3d outputs for the domain point's embedding in $\mathbb{R}^3$.
The architecture uses an initial spectral action to intrinsically enforce domain identifications, and a balance of losses to enforce regularity and smoothness of the surface (including domain gluing for disconnected constructions).

The training process minimises a Willmore loss, as a Monte-Carlo estimation of the Willmore energy over each training batch, mimicking the Willmore flow as a neural training process.
The Willmore energy is computed via the mean curvature, using \texttt{pytorch}'s inbuilt autograd functionality to efficiently compute partial derivatives at millions of points, treating the Willmore energy as a PINN-style loss.

This novel neural Willmore flow method is tested on surfaces of genus $\{0, 1, 2\}$, where the former genus 0 case validates the proven results in the literature for the round sphere.
Following this, application to the famous recently proven result \cite{MarquesNeves2014} for genus 1, surprisingly manages to numerically reproduce the convergence to the Clifford torus minima, for varying $\tau$ values, and with supporting visualisations. 
Finally, the latter case of genus 2 is an open problem, and this approach successfully produces a surface with Willmore energy substantially lower than the naive catenoid bridge connection between two Clifford tori, yet not as low as the currently conjectured $\xi_{2,1}$ Lawson surface minimum.
The smooth decreases in Willmore energy motivate further testing and refined training towards this Willmore minima.
Irrespectively, this method is cemented as an insightful and novel way to search for Willmore minimal surfaces using tools from the modern AI era.

Excitingly, this work sets up an array of potentially impactful developments.
As well as refining the genus 2 learning in search for a true lower Willmore minima, the code is well positioned for immediate extension to higher genuses via this tori gluing construction.
Whilst for closed connected orientable 2d surfaces the diffeomorphism classes are defined by genus alone, in higher dimensions the classification becomes far more technical, and in each case (assuming a known manifold can be used for initialisation), this codebase can be simply extended for application there.
Additionally, beyond just embeddings in $\mathbb{R}^3$ (or $\mathbb{R}^n$), one may consider embedding in other space forms of either positive $\mathbb{S}^n$ or negative $\mathbb{H}^n$ curvature.
The architectures themselves also have no obstruction to application on more general immersions, which then permits investigations into more complicated surfaces with self-intersection, such as the Klein bottle.
For any of these suggestions, the authors welcome readers to develop the ideas using this work's code in their own collaborations.
Finally, whilst the central novel loss in this work is the Willmore energy, its construction via the mean curvature opens up potential for building losses associated to constant mean curvature, which the authors are presently pursuing. 

\section*{Acknowledgements}
The authors wish to thank Lucca Delboni, and Jakob Stein, for useful comments.\\
\textbf{EH} is supported by São Paulo Research Foundation (FAPESP) grant 2024/18994-7. \\
\textbf{HSE} is supported by grants: FAPESP 2020/09838-0, with BI0S: Brazilian Institute of Data Science, 
FAPESP 2021/04065-6, with BRIDGES: Brazil-France interplays in Gauge Theory, extremal structures and stability,
FAPESP 2024/00923-6, with CBG: Brazilian Centre for Geometry. 307145/2025-5 level PQ-A, with the Brazilian National Council for Scientific and Technological Development (CNPq).\\
\textbf{TS} is supported by FAPESP grant 2022/09891-4.

\section*{Data Availability}
The code used in this work is available at the respective repository: \url{https://github.com/edhirst/WillmorePINN}. 

\appendix
\section{A brief overview of machine learning}
\label{appendix:ml}

This appendix gives a non-technical overview of machine learning for the
reader encountering these ideas for the first time. The precise formulation
used in this work is developed in \S\ref{sec:pinns_bkg}; here we aim only to
convey the main ideas and fix basic terminology. For comprehensive treatments
see~\cite{Bishop2016, Goodfellow-et-al-2016, prince2023understanding}.

Machine learning is, broadly, the study of algorithms that improve their
performance on a task by leveraging data. Rather than writing down an explicit
procedure to solve a problem, one specifies a flexible family of candidate
solutions and an objective that measures how well each candidate performs; an
automated search then selects the best candidate. From a mathematician's
perspective, this is simply an optimisation problem, distinguished mainly by
the scale of the parameter space and by the use of statistical, data-driven
objectives.

The candidate solutions are called \emph{models}. The dominant model class
today is the \emph{neural network}: a parametrised function built by
alternating affine maps with simple fixed nonlinearities (such as
$x \mapsto \max(0,x)$). Despite this elementary construction, neural networks
are \emph{universal approximators} --- given enough parameters, they can
approximate any continuous function on a compact domain to arbitrary
precision~\cite{Cybenko1989, Hornik1991}. The parameters of a neural network
(often numbering in the thousands to millions) are called \emph{weights}, and
the process of selecting them is called \emph{training}.

Training is driven by a \emph{loss function} $\mathcal{L}(\theta)$, a scalar
objective defined on the weight space that quantifies how far the current model is from the desired behaviour. The optimisation is performed by \emph{gradient descent}: one computes $\nabla_\theta \mathcal{L}$ and takes a small step in the direction of steepest decrease. This gradient is obtained efficiently by \emph{automatic differentiation} (also called \emph{backpropagation}), an algorithmic application of the chain rule that is exact up to floating-point precision. Modern software libraries such as \texttt{pytorch}~\cite{torch} make this routine essentially automatic, even for models with millions of parameters.

The nature of the loss function determines the \emph{learning paradigm}. In
\emph{supervised learning}, the most classical setting, one possesses a
dataset of input--output pairs and the loss measures discrepancy between the
model's predictions and the known outputs. In \emph{unsupervised learning},
no target outputs are provided and the loss instead captures some intrinsic
structure of the data (e.g.\ clustering or compression).

A third paradigm, and the one most relevant to this work, replaces the
data-fitting loss with a \emph{physics-informed} objective: instead of
comparing the model to pre-computed answers, the loss encodes the governing
equations or variational principle of the problem directly. This is the idea
behind \emph{physics-informed neural networks}
(PINNs)~\cite{Raissi2019}. Because the loss involves derivatives of the
model -- computed exactly by automatic differentiation -- and because sample
points in the domain can be generated freely, a PINN requires no external
dataset at all: the ``data'' is the geometry of the problem itself. Further
technical details of this approach, and its application to the Willmore
energy, are discussed in \S\ref{sec:pinns_bkg}.

A first time coder can get started using this codebase, without need to set up environments or GitHub, by using the accompanying \href{https://mybinder.org/v2/gh/edhirst/WillmorePINN/HEAD?urlpath=%2Fdoc%2Ftree%2Fdemo.ipynb}{Binder notebook}.
Here, one can alter hyperparameters that change the training objective and run directly in the cloud without any local set up.

\section{Twisted genus 1}\label{app:g1_twisted}
This appendix looks to demonstrate the greater generality of the accompanying codebase, which allows for general defining parameter values across the varying genuses considered.

Here, an example of the neural Willmore flow is shown for a genus 1 initialisation surface with $\tau = 0.1 + 0.1i$.
With this value the fundamental domain has shearing\footnote{Note the fundamental domain is sampled and coloured in the same way, then scaled by the $\tau$ factor to produce the true sample used as architecture input.}, leading to a twisting in the embedded torus.
Results for an equivalent run (with the same training hyperparameters) to that performed in §\ref{sec:results_g1} are given below.

\begin{figure}[H]
    \centering
    \begin{subfigure}[t]{0.32\textwidth}
        \centering
        \includegraphics[width=\linewidth]{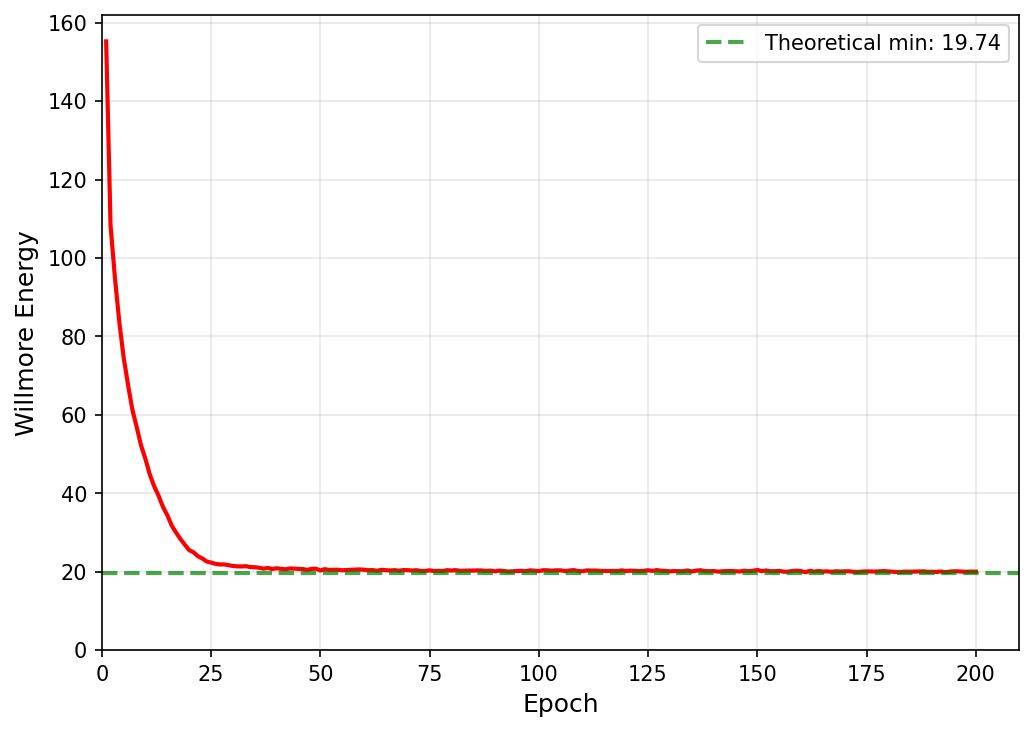}
        \caption{Willmore Loss}
        \label{fig:g1_twist_willmore}
    \end{subfigure}
    \hfill
    \begin{subfigure}[t]{0.32\textwidth}
        \centering
        \includegraphics[width=\linewidth]{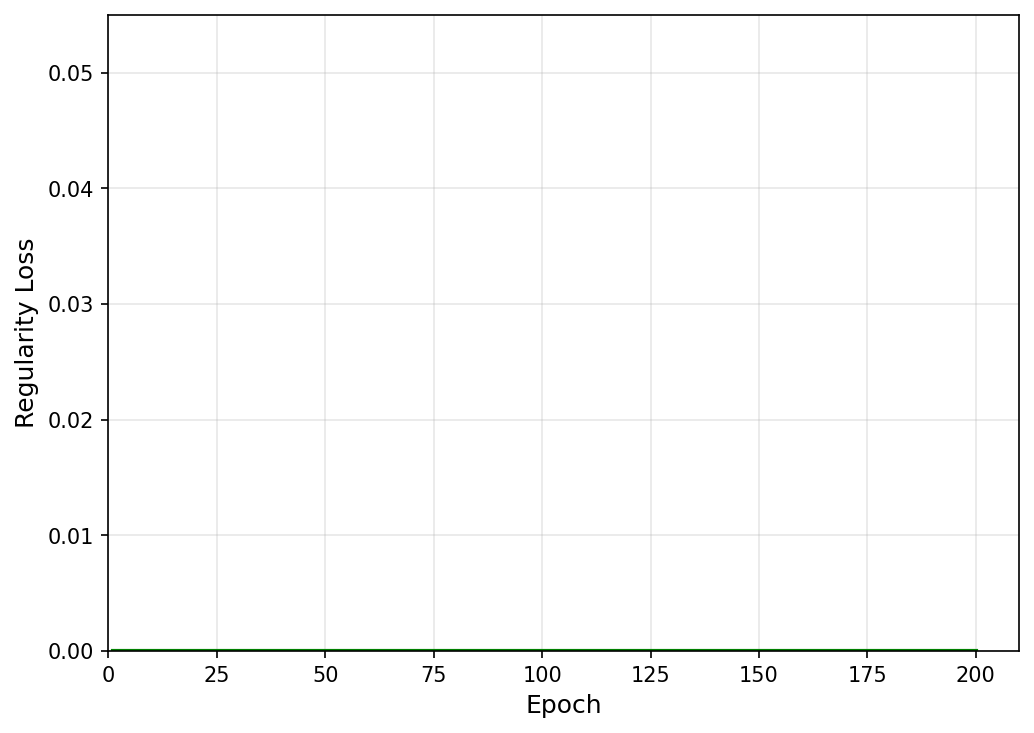}
        \caption{Regularity Loss}
        \label{fig:g1_twist_regularity}
    \end{subfigure}
    \hfill
    \begin{subfigure}[t]{0.32\textwidth}
        \centering
        \includegraphics[width=\linewidth]{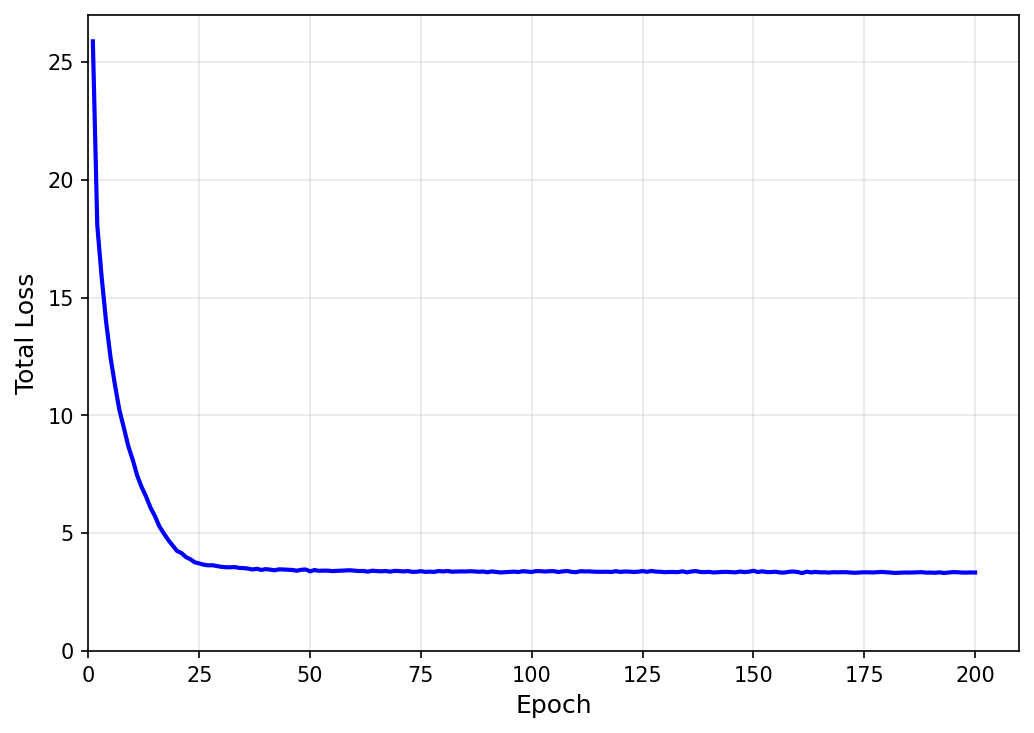}
        \caption{Total Loss}
        \label{fig:g1_twist_total}
    \end{subfigure}
    \caption{Training losses for the twisted genus 1 setup. Willmore loss reaches the expected theoretical minimum for the Clifford torus from the $\tau = 0.1 + 0.1i$ twisted start point; regularity loss stays at exactly 0 throughout (below the ReLU-gated thresholds) affirming a smooth consistent surface.}
    \label{fig:g1twist_losses}
\end{figure}

\begin{figure}[H]
    \centering
    \begin{subfigure}[t]{0.32\textwidth}
        \centering
        \includegraphics[width=0.8\linewidth]{figures/genus1/fundamental_domain.png}
        \caption{Fundamental Domain}
        \label{fig:g1t_fd}
    \end{subfigure}
    \hfill
    \begin{subfigure}[t]{0.32\textwidth}
        \centering
        \includegraphics[width=\linewidth]{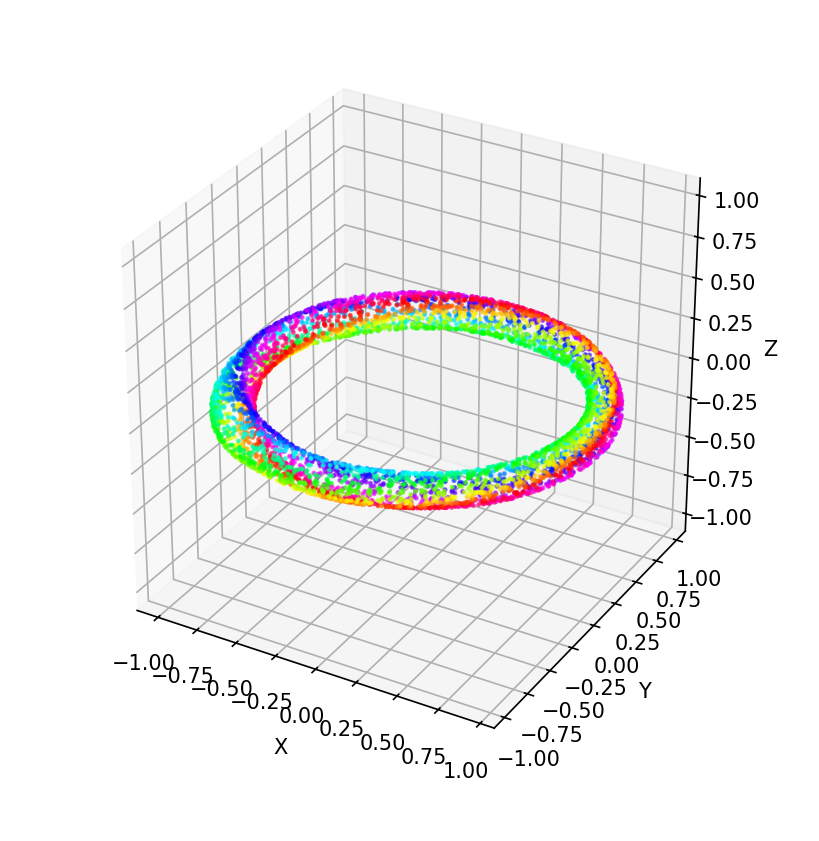}
        \caption{Epoch 0, $\widehat{\mathcal{W}} = 275.52$}
        \label{fig:g1t_e0}
    \end{subfigure}
    \hfill
    \begin{subfigure}[t]{0.32\textwidth}
        \centering
        \includegraphics[width=\linewidth]{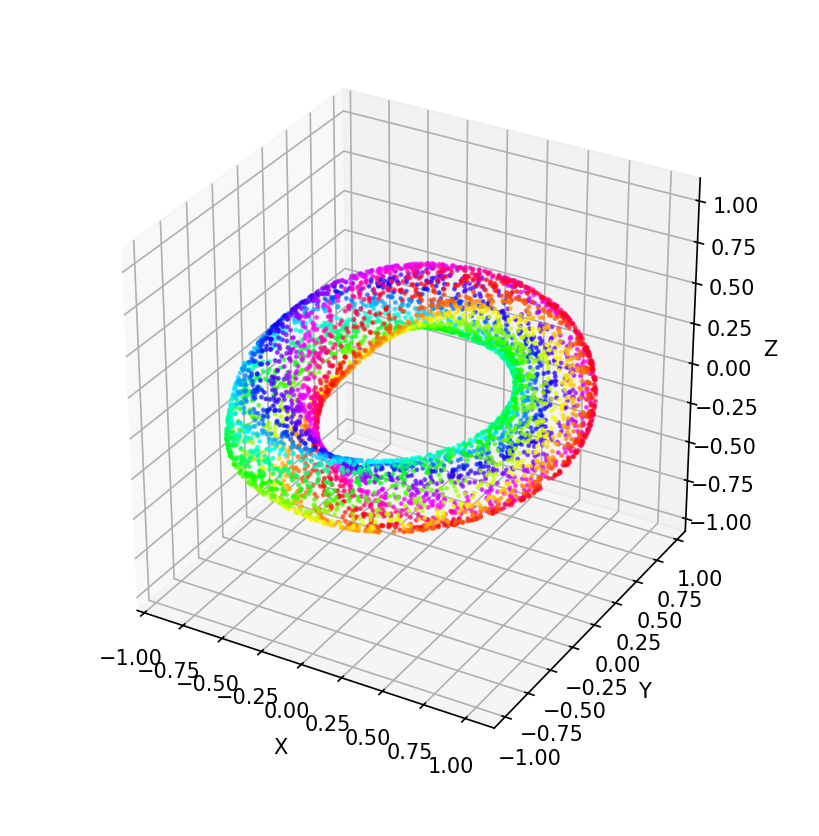}
        \caption{Epoch 10, $\widehat{\mathcal{W}} = 48.68$}
        \label{fig:g1t_e10}
    \end{subfigure}\\
    \begin{subfigure}[t]{0.32\textwidth}
        \centering
        \includegraphics[width=\linewidth]{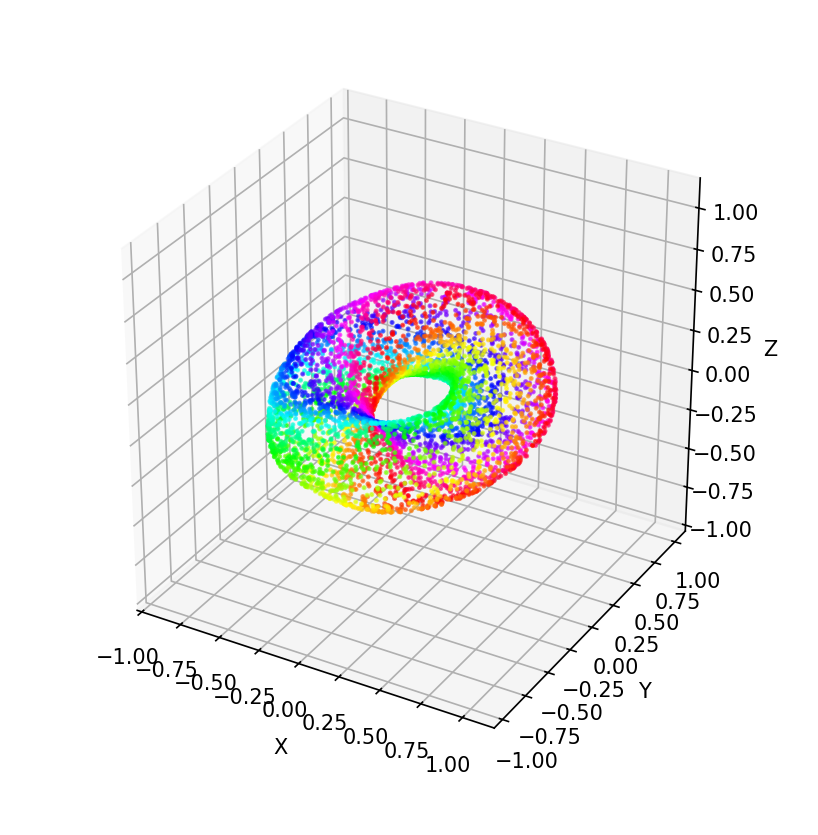}
        \caption{Epoch 20, $\widehat{\mathcal{W}} = 25.48$}
        \label{fig:g1t_e20}
    \end{subfigure}
    \hfill
    \begin{subfigure}[t]{0.32\textwidth}
        \centering
        \includegraphics[width=\linewidth]{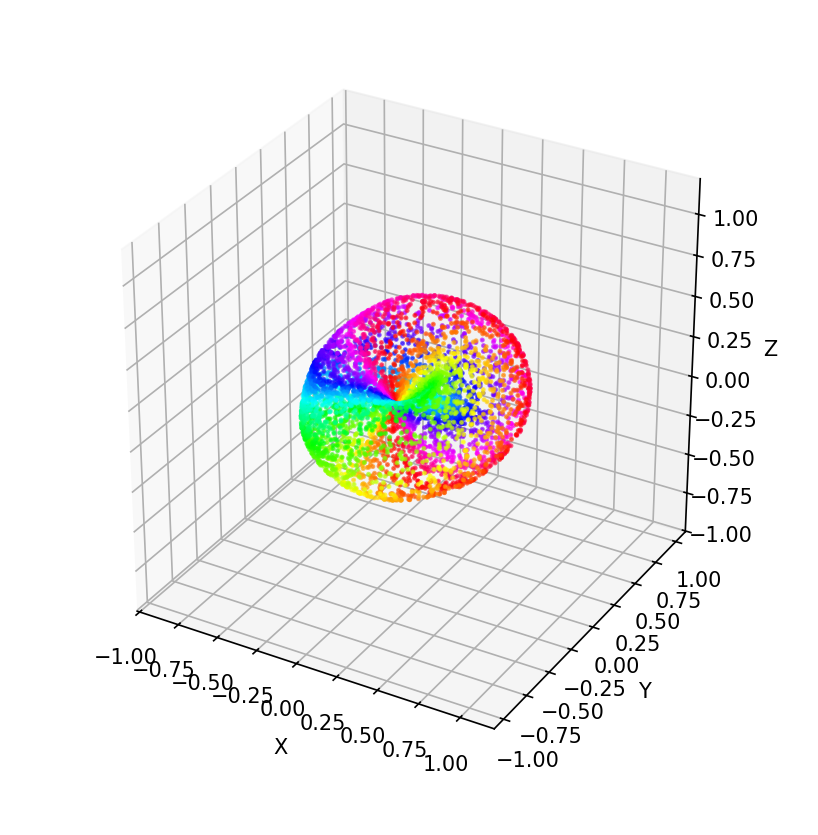}
        \caption{Epoch 100, $\widehat{\mathcal{W}} = 20.14$}
        \label{fig:g1t_e100}
    \end{subfigure}
    \hfill
    \begin{subfigure}[t]{0.32\textwidth}
        \centering
        \includegraphics[width=\linewidth]{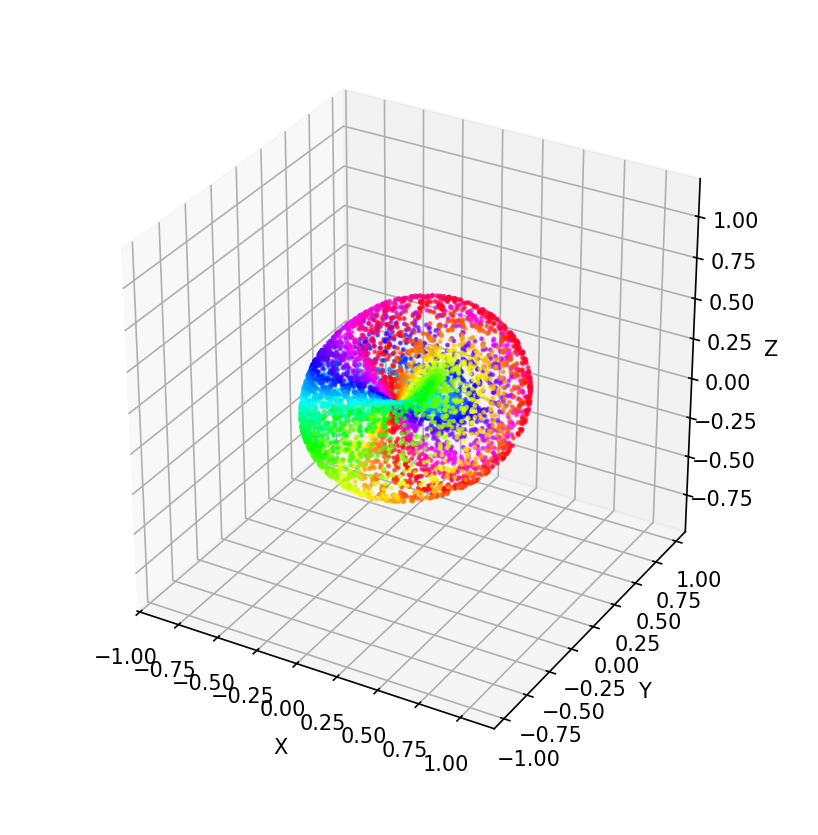}
        \caption{Epoch 200, $\widehat{\mathcal{W}} = 20.00$}
        \label{fig:g1t_e200}
    \end{subfigure}
    \caption{Visualisation of the embedding evolution over training for the genus 1 setup. The fundamental domain is shown in (a) where the $(u,v)$ coordinates denote the architecture input (after conversion to Fourier modes), then (b)-(f) are embedding outputs of the model at various epochs through training from $\tau = 0.1 + 0.1i$, using point colours matching the fundamental domain input, with values of the Willmore energy $\widehat{\mathcal{W}}$.}
    \label{fig:g1twist_visuals}
\end{figure}

Performance is equivalently strong, were losses in Figure \ref{fig:g1twist_losses} show approaching of the minimal Clifford torus, yet with a twist as shown in the embedding visualisations of Figure \ref{fig:g1twist_visuals}. 
It is worth noting that all testing with varied $\tau$ values showed similar behaviour, and equivalently for other genuses.

\addcontentsline{toc}{section}{References}
\bibliographystyle{utphys}
\bibliography{references}{}

\providecommand{\href}[2]{#2}\begingroup\raggedright\begin{thebibliography}{10}

\bibitem{Willmore1965}
T.~J. Willmore, ``Note on embedded surfaces,'' {\em An. {\c S}tiin{\c t}. Univ. ``Al. I. Cuza'' Ia{\c s}i Sec{\c t}. I a Mat. (N.S.)} {\bfseries 11B} (1965) 493--496.

\bibitem{MarquesNeves2014}
F.~C. Marques and A.~Neves, ``Min-max theory and the willmore conjecture,'' {\em Annals of Mathematics} {\bfseries 179} no.~2, (2014) 683--782.

\bibitem{Hirst:2025seh}
E.~Hirst, T.~S. Gherardini, and A.~G. Stapleton, ``{AInstein: Numerical Einstein Metrics via Machine Learning},'' \href{http://arxiv.org/abs/2502.13043}{{\ttfamily arXiv:2502.13043 [hep-th]}}.

\bibitem{Cortes:2026kfx}
G.~Cort{\'e}s, M.~Esteban-Casadevall, Y.~Feng, J.~Henkel, E.~Hirst, T.~S. Gherardini, and A.~G. Stapleton, ``{A Machine Learning Approach to the Nirenberg Problem},'' \href{http://arxiv.org/abs/2602.12368}{{\ttfamily arXiv:2602.12368 [cs.LG]}}.

\bibitem{douglas2020numerical}
D.~Michael, L.~Subramanian, and Q.~Yidi, ``{Numerical Calabi--Yau metrics from holomorphic networks},'' {\em Proceedings of Machine Learning Research} {\bfseries 145} (2020) 223--252, \href{http://arxiv.org/abs/2012.04797}{{\ttfamily arXiv:2012.04797 [hep-th]}}.

\bibitem{Larfors:2021pbb}
M.~Larfors, A.~Lukas, F.~Ruehle, and R.~Schneider, ``{Learning Size and Shape of Calabi-Yau Spaces},'' 11, 2021.

\bibitem{Gerdes:2022nzr}
M.~Gerdes and S.~Krippendorf, ``{{CYJAX}: A package for {Calabi--Yau} metrics with {JAX}},'' \href{http://dx.doi.org/10.1088/2632-2153/acdc84}{{\em Mach. Learn. Sci. Tech.} {\bfseries 4} no.~2, (2023) 025031}, \href{http://arxiv.org/abs/2211.12520}{{\ttfamily arXiv:2211.12520 [hep-th]}}.

\bibitem{Berglund:2022gvm}
P.~Berglund, G.~Butbaia, T.~H\"uubsch, V.~Jejjala, D.~Mayorga Pe\~na, C.~Mishra, and J.~Tan, ``{Machine-learned {Calabi--Yau} metrics and curvature},'' \href{http://dx.doi.org/10.4310/ATMP.2023.v27.n4.a3}{{\em Adv. Theor. Math. Phys.} {\bfseries 27} no.~4, (2023) 1107--1158}, \href{http://arxiv.org/abs/2211.09801}{{\ttfamily arXiv:2211.09801 [hep-th]}}.

\bibitem{douglas2024harmonic}
M.~R. Douglas, D.~Platt, and Y.~Qi, ``{Harmonic 1-forms on real loci of Calabi--Yau manifolds},'' May, 2024.
\newblock \url{https://arxiv.org/abs/2405.19402}.

\bibitem{Heyes:2026rch}
E.~Heyes, E.~Hirst, H.~N.~S. Earp, and T.~S.~R. Silva, ``{Neural and numerical methods for $\mathrm{G}_2$-structures on contact Calabi-Yau 7-manifolds},'' \href{http://arxiv.org/abs/2602.12438}{{\ttfamily arXiv:2602.12438 [math.DG]}}.

\bibitem{GomezSerranoSurvey}
J.~G{\'o}mez-Serrano, ``Computer-assisted proofs in pde: A survey,'' {\em Notices of the AMS} {\bfseries 66} no.~3, (2019) 298--310.

\bibitem{wang2023}
Y.~Wang, C.-Y. Lai, J.~G\'omez-Serrano, and T.~Buckmaster, ``Asymptotic self-similar blow-up profile for three-dimensional axisymmetric {Euler} equations using neural networks,'' \href{http://dx.doi.org/10.1103/PhysRevLett.130.244002}{{\em Phys. Rev. Lett.} {\bfseries 130} (Jun, 2023) 244002}. \url{https://link.aps.org/doi/10.1103/PhysRevLett.130.244002}.

\bibitem{platt_nirenberg}
D.~Platt, ``Non-uniqueness and symmetries for the nirenberg problem using computer assistance,'' 2026.
\newblock \url{https://arxiv.org/abs/2603.29544}.

\bibitem{ZhouYe2023MinimalSurfacePINN}
S.~Zhou and X.~Ye, ``Approximating high-dimensional minimal surfaces with physics-informed neural networks,'' 2023.

\bibitem{Hashimoto:2025zmi}
K.~Hashimoto, K.~Kyo, M.~Murata, G.~Ogiwara, and N.~Tanahashi, ``{Physics-informed neural network solves minimal surfaces in curved spacetime},'' \href{http://dx.doi.org/10.1088/2632-2153/ae3050}{{\em Mach. Learn. Sci. Tech.} {\bfseries 7} no.~1, (2026) 015013}, \href{http://arxiv.org/abs/2509.10866}{{\ttfamily arXiv:2509.10866 [hep-th]}}.

\bibitem{torch}
A.~Paszke, S.~Gross, F.~Massa, A.~Lerer, J.~Bradbury, G.~Chanan, T.~Killeen, Z.~Lin, N.~Gimelshein, L.~Antiga, A.~Desmaison, A.~Köpf, E.~Yang, Z.~DeVito, M.~Raison, A.~Tejani, S.~Chilamkurthy, B.~Steiner, L.~Fang, J.~Bai, and S.~Chintala, ``Pytorch: An imperative style, high-performance deep learning library,'' 2019.
\newblock \url{https://arxiv.org/abs/1912.01703}.

\bibitem{Thomsen1923}
G.~Thomsen, ``{\"U}ber konforme {G}eometrie {I}: {G}rundlagen der konformen {F}l\"achentheorie,'' {\em Abh.\ Math.\ Semin.\ Univ.\ Hambg.} {\bfseries 3} (1923) 31--56.

\bibitem{Blaschke1929}
W.~Blaschke, {\em Vorlesungen \"uber {D}ifferentialgeometrie, {I}{I}{I}}.
\newblock Springer, Berlin, 1929.

\bibitem{White1973}
J.~H. White, ``A global invariant of conformal mappings in space,'' {\em Proc.\ Amer.\ Math.\ Soc.} {\bfseries 38} no.~1, (1973) 162--164.

\bibitem{Halverson:2023ndu}
J.~Halverson and F.~Ruehle, ``{Metric flows with neural networks},'' \href{http://dx.doi.org/10.1088/2632-2153/ad8533}{{\em Mach. Learn. Sci. Tech.} {\bfseries 5} no.~4, (2024) 045020}, \href{http://arxiv.org/abs/2310.19870}{{\ttfamily arXiv:2310.19870 [hep-th]}}.

\bibitem{LiYau1982}
P.~Li and S.-T. Yau, ``A new conformal invariant and its applications to the {W}illmore conjecture and the first eigenvalue of compact surfaces,'' {\em Invent.\ Math.} {\bfseries 69} no.~2, (1982) 269--291.

\bibitem{Simon1993}
L.~Simon, ``Existence of surfaces minimizing the willmore functional,'' {\em Communications in Analysis and Geometry} {\bfseries 1} no.~2, (1993) 281--326.

\bibitem{BauerKuwert2003}
M.~Bauer and E.~Kuwert, ``Existence of minimizing {W}illmore surfaces of prescribed genus,'' {\em Int.\ Math.\ Res.\ Not.} no.~10, (2003) 553--576.

\bibitem{Riviere2008}
T.~Rivi\`ere, ``Analysis aspects of {W}illmore surfaces,'' {\em Invent.\ Math.} {\bfseries 174} no.~1, (2008) 1--45.

\bibitem{Kusner1989}
R.~Kusner, ``Comparison surfaces for the {W}illmore problem,'' {\em Pacific J.\ Math.} {\bfseries 138} no.~2, (1989) 317--345.

\bibitem{Lawson1970}
H.~B. Lawson, ``Complete minimal surfaces in {$S^3$},'' {\em Ann.\ of Math.\ (2)} {\bfseries 92} no.~3, (1970) 335--374.

\bibitem{Hsu1992}
L.~Hsu, R.~Kusner, and J.~Sullivan, ``Minimizing the squared mean curvature integral for surfaces in space forms,'' {\em Experiment.\ Math.} {\bfseries 1} no.~3, (1992) 191--207.

\bibitem{Goodfellow2016}
I.~Goodfellow, Y.~Bengio, and A.~Courville, {\em Deep Learning}.
\newblock MIT Press, 2016.

\bibitem{Rumelhart1986}
D.~E. Rumelhart, G.~E. Hinton, and R.~J. Williams, ``Learning representations by back-propagating errors,'' {\em Nature} {\bfseries 323} (1986) 533--536.

\bibitem{Kingma2015}
D.~P. Kingma and J.~Ba, ``Adam: A method for stochastic optimization,'' \href{http://arxiv.org/abs/1412.6980}{{\ttfamily arXiv:1412.6980 [cs.LG]}}. \url{https://arxiv.org/abs/1412.6980}.

\bibitem{Cybenko1989}
G.~Cybenko, ``Approximation by superpositions of a sigmoidal function,'' {\em Math.\ Control Signals Systems} {\bfseries 2} no.~4, (1989) 303--314.

\bibitem{Hornik1991}
K.~Hornik, ``Approximation capabilities of multilayer feedforward networks,'' {\em Neural Netw.} {\bfseries 4} no.~2, (1991) 251--257.

\bibitem{Barron1993}
A.~R. Barron, ``Universal approximation bounds for superpositions of a sigmoidal function,'' {\em IEEE Trans.\ Inform.\ Theory} {\bfseries 39} no.~3, (1993) 930--945.

\bibitem{DeVore2021}
R.~DeVore, B.~Hanin, and G.~Petrova, ``Neural network approximation,'' {\em Acta Numer.} {\bfseries 30} (2021) 327--444.

\bibitem{Yarotsky2017}
D.~Yarotsky, ``Error bounds for approximations with deep {R}e{LU} networks,'' {\em Neural Netw.} {\bfseries 94} (2017) 103--114.

\bibitem{Raissi2019}
M.~Raissi, P.~Perdikaris, and G.~E. Karniadakis, ``Physics-informed neural networks: {A} deep learning framework for solving forward and inverse problems involving nonlinear partial differential equations,'' {\em J.\ Comput.\ Phys.} {\bfseries 378} (2019) 686--707.

\bibitem{Karniadakis2021}
G.~E. Karniadakis, I.~G. Kevrekidis, L.~Lu, P.~Perdikaris, S.~Wang, and L.~Yang, ``Physics-informed machine learning,'' {\em Nat.\ Rev.\ Phys.} {\bfseries 3} no.~6, (2021) 422--440.

\bibitem{Baydin2018}
A.~G. Baydin, B.~A. Pearlmutter, A.~A. Radul, and J.~M. Siskind, ``Automatic differentiation in machine learning: a survey,'' {\em J.\ Mach.\ Learn.\ Res.} {\bfseries 18} no.~153, (2018) 1--43.

\bibitem{Dziuk2008}
G.~Dziuk and C.~M. Elliott, ``Finite elements on evolving surfaces,'' {\em IMA J.\ Numer.\ Anal.} {\bfseries 27} no.~2, (2007) 262--292.

\bibitem{Crane2013}
K.~Crane, U.~Pinkall, and P.~Schr\"{o}der, ``Robust fairing via conformal curvature flow,'' {\em ACM Trans.\ Graph.} {\bfseries 32} no.~4, (2013) 61:1--61:10.

\bibitem{lawson_visuals}
{GeometrieWerkstatt}, ``Lawson surfaces.'' \url{http://geometriewerkstatt.com/LawsonSurfaces.html}.

\bibitem{Bishop2016}
C.~M.~Bishop, {\em Pattern Recognition and Machine Learning}.
\newblock Information Science and Statistics. Springer, New York, NY, Aug., 2016.

\bibitem{Goodfellow-et-al-2016}
I.~Goodfellow, Y.~Bengio, and A.~Courville, {\em Deep Learning}.
\newblock MIT Press, Cambridge, MA, 2016.
\newblock \url{http://www.deeplearningbook.org}.

\bibitem{prince2023understanding}
S.~J. Prince, {\em Understanding Deep Learning}.
\newblock The MIT Press, Cambridge, MA, 2023.
\newblock \url{http://udlbook.com}.

\end{thebibliography}\endgroup

\end{document}